\magnification1000
 \def\oui{oui}
 \ifx\textures\oui
 \input CrayolaColors 
 \long\def\rge#1{\Red#1\Black}

 \else

  
  %
  %

  %
  %
  \let\footnotea=\footnote
  \def\anote#1#2{\footnotea{\hbox{$^{#1}$}}{\eightpoint#2}}  
  \catcode`@=12 

 \def\defrefnote#1{\definexref{#1}{{\the\footnotenumber}}{refnotes}}

  %
  %


 \input eplain.tex
\makeatletter
\def\numberedfootnote{%
ÊÊ\global\advance\footnotenumber by 1
ÊÊ\@eplainfootnote{{\number\footnotenumber}}%
}%
\def\makecolumns#1/#2 {\par \begingroup
ÊÊ \@columndepth = #1
ÊÊ \advance\@columndepth by -1
ÊÊ \divide \@columndepth by #2
ÊÊ \advance\@columndepth by 1
ÊÊ \@linestogoincolumn = \@columndepth
ÊÊ \@linestogo = #1
ÊÊ \currentcolumn = 1
ÊÊ \def\@endcolumnactions{%
ÊÊÊÊÊÊ\ifnum \@linestogo<2
ÊÊÊÊÊÊÊÊ \the\crtok \egroup \endgroup \par 
ÊÊÊÊÊÊ\else
ÊÊÊÊÊÊÊÊ \global\advance\@linestogo by -1
ÊÊÊÊÊÊÊÊ \ifnum\@linestogoincolumn<2
ÊÊÊÊÊÊÊÊÊÊÊÊ\global\advance\currentcolumn by 1
ÊÊÊÊÊÊÊÊÊÊÊÊ\global\@linestogoincolumn = \@columndepth
ÊÊÊÊÊÊÊÊÊÊÊÊ\the\crtok
ÊÊÊÊÊÊÊÊ \else
ÊÊÊÊÊÊÊÊÊÊÊÊ&\global\advance\@linestogoincolumn by -1
ÊÊÊÊÊÊÊÊ \fi
ÊÊÊÊÊÊ\fi
ÊÊ }%
ÊÊ \makeactive\^^M
ÊÊ \letreturn \@endcolumnactions
ÊÊ \@columnwidth = \hsize
ÊÊÊÊ \advance\@columnwidth by -\parindent
ÊÊÊÊ \divide\@columnwidth by #2
ÊÊ \penalty\abovecolumnspenalty
ÊÊ \noindent 
ÊÊ \valign\bgroup
ÊÊÊÊ &\hbox to \@columnwidth{\strut \hsize = \@columnwidth ##\hfil}\cr
}%
\makeatother

\lefteqnumbers
   \def\testd{oui}
   \def\choixlat{\ifx\numadroite\testd\righteqnumbers
            \else  \lefteqnumbers\fi}
    \choixlat

\catcode`@=\letter
\def\@eplainfootnote#1{\let\@sf\empty 
  \ifhmode\edef\@sf{\spacefactor\the\spacefactor}\/\fi
  \global\advance\hlfootlabelnumber by 1
  \hlstart@impl{foot}{\hlfootlabel}%
  \hldest@impl{footback}{\hlfootbacklabel}%
  \hbox{$^{(#1)}$}%
  \hlend@impl{foot}%
  \@sf\vfootnote{#1.}%
}%
\catcode`@=\other

  \interfootnoteskip=0pt
  \let\note=\numberedfootnote
  \everyfootnote={\eightpoint\leftskip=5truemm\rightskip5truemm}
  
  \hsize150truemm\vsize 240truemm\hoffset=5truemm
  \def\dimstand{\hsize 150truemm\vsize 240truemm\hoffset=5truemm\voffset=0truemm}
  
  \pretolerance=500\tolerance=1000\brokenpenalty=5000
  \parindent3mm
  
  \countdef\temps=170
  \temps=\time
  \countdef\nminutes=171{\nminutes = \time}
  \countdef\nheure=172
  \def\heure{\begingroup                     
     \temps = \time \divide\temps by 60
     \nheure = \temps                        
     \nminutes = \time
     \multiply\temps by 60
     \advance\nminutes by -\temps            
     \ifnum\nminutes<10 \toks1 = {0}%
     \else\toks1 = {}%
     \fi
     \number\nheure h\the\toks1 \number\nminutes  
  \endgroup}%

  \newcount\chstart
  \chstart=\pageno
 \headline={\ifnum\pageno=\chstart {\hfill} \else {\hss \tenrm --\ \folio\ --\hss}\fi}
  \footline={\hfill}
  \normalbaselines
  \frenchspacing
    \def\dater{\vglue-10mm\rightline{(\the\day/\the\month/\the\year)}}
  \def\dateheure{\vglue-10mm\rightline{(\the\day/\the\month/\the\year,\ \heure)}}

  \newif\ifpagetitre \pagetitretrue
\newtoks\hautpagetitre \hautpagetitre={\hfill}
\newtoks\baspagetitre \baspagetitre={\hfill}
\newtoks\auteurcourant \auteurcourant={\hfill}
\newtoks\titrecourant \titrecourant={\hfill}
\newtoks\hautpagegauche
\newtoks\hautpagedroite
\newtoks\hautpagemilieu
\hautpagemilieu={\tenrm\hfil -- \folio\ -- \hfil}
\hautpagegauche={\ifx\midfolio\oui\the\hautpagemilieu\else\tenrm\folio\hfill\the\auteurcourant\hfill\fi}
\hautpagedroite={\ifx\midfolio\oui\the\hautpagemilieu\else\hfill\the\titrecourant\hfill\tenrm\folio\fi}
\newtoks\baspagegauche \baspagegauche={\hfil}
\newtoks\baspagedroite \baspagedroite={\hfil}
\headline={\ifpagetitre\the\hautpagetitre
\else\ifodd\pageno\the\hautpagedroite\else\the\hautpagegauche\fi\fi }
\footline={\ifpagetitre\the\baspagetitre
\else\ifodd\pageno\the\baspagedroite
\else\the\baspagegauche\fi\fi \global\pagetitrefalse}

\def\pageblanche{\vfill\eject\pagetitretrue
\null\vfill\eject
\pagetitretrue
}
\def\chgtpage{\ifodd\pageno \else
\pageblanche \fi \pagetitretrue\titreun=0\footnotenumber=0}

\def\chgtpageincrtitreun{\ifodd\pageno \else
\pageblanche \fi \pagetitretrue\footnotenumber=0}

\def\majnombres{\ifodd\pageno \else
\pageblanche \fi \pagetitretrue\hautpoly\titreun=0\footnotenumber=0}

\def\hautspages#1#2{\auteurcourant={\ninepcap#1}\titrecourant={\nineit#2}}

\ifnum\chstart=\pageno \pagetitretrue\fi
  

  \def\PAR{\par}
  
  \def\leftnote#1{\vadjust{\setbox1=\vtop{\hsize 20mm\parindent=0pt\eightpoint
  \baselineskip=9pt\rightskip=4mm plus 4mm\vskip-4.7mm#1}\hbox{\kern-2cm\smash{\box1}}}}

  
  \def\raggedcenter{\leftskip=20pt plus 10em  
       \rightskip=\leftskip 
        \parfillskip=0pt 
         \spaceskip=.3333em \xspaceskip=.5em 
          \pretolerance=9999 \tolerance=9999
           \hyphenpenalty=9999 \exhyphenpenalty=9999 }
           
  \def\titrecentre#1{{\parindent0mm\raggedcenter
       \spaceskip=.6em plus .2em minus .2em\xspaceskip=.6em plus .2em minus .2em
        \tit#1\par}}
        


  \def\oui{oui}
  
   \def\fontetitreun{\ifx\paradouze\oui\douzepts\gpdouze\twelvebf\textfont1=\twelveib\else
\quatorzepts\gpquatorze\fourteenbf\fi}

\def\fontetitreunl{\douzepts\textfont1=\twelveib\scriptfont1=\tenib\fourteenti}
 
 \def\fontetitredeux{\textfont1=\eleveni\ifx\paradouze\oui\onzepts\scriptfont1=\ninei\elevenit\else
                        \douzepts\twelveit\fi}
 
   \def\fontetitredeuxb{\ifx\paradouze\oui\onzepts\eleventi\gponze\textfont1=\elevenib\scriptfont1=\nineib
                         \else\douzepts\twelveti\scriptfont1=\twelveib\scriptfont1=\tenib\gpdouze\fi}
                         
\def\fontetitredeuxl{\onzepts\textfont1=\elevenbf\scriptfont1=\ninebf\twelvebf}
  
\def\fontetitretrois{\textfont0=\elevenrm\scriptfont0=\eightrm\textfont1=\eleveni
                      \scriptfont1=\eighti\scriptscriptfont1=\sixi\elevenit}
                      
\def\fontetitrequatre{\textfont0=\elevenrm\scriptfont0=\eightrm\textfont1=\eleveni
                      \scriptfont1=\eighti\scriptscriptfont1=\sixi\elevenrm}
  
  \newcount\titreun\titreun=0
  \newcount\titredeux\titredeux=0
  \newcount\titretrois\titretrois=0
  \newcount\titrequatre\titrequatre=0
  \newcount\enonce\enonce=0
  
  \def\incr#1{\global\advance#1 by 1 {\the #1}}
  \def\avance#1{\global\advance#1 by 1}
  \def\init#1{\global#1=0}
  
  \long\def\Indentation#1#2{\setbox10=\hbox{\fontetitreun#1}
  	                    \ifdim\wd10 < 4mm
                         \setbox10=\hbox to 4mm{\box10\hfill}
                       \else\ifdim\wd10 < 6mm
                         \setbox10=\hbox to 6mm{\box10\hfill}
  	                    \else\ifdim\wd10 < 8mm
                         \setbox10=\hbox to 8mm{\box10\hfill}
                       \else\ifdim\wd10 < 12mm
                         \setbox10=\hbox to 12mm{\box10\hfill}
                       \fi\fi\fi\fi
                       \dimen10=\hsize
                       \advance \dimen10 by -\wd10
                       \noindent \box10 %
                       \ignorespaces
                       \hbox{\vtop{\hsize=\dimen10\raggedright\noindent\fontetitreun#2}}}

  \long\def\paraun#1{\removelastskip\par\medskip\goodbreak\vskip0pt plus.01\vsize\penalty-100
                \vskip0pt plus-.01\vsize
  	              \init{\titredeux}\ifnum\optionparag=1{\init\eqnumber\init\enonce}\else{}\fi
                  \goodbreak{\fontetitreun
  	                \Indentation{\incr{\titreun}.\ }{\fontetitreun #1\par}}\nobreak\medskip}

 %
 %
 \long\def\paraunc#1{\removelastskip\par\bigskip\goodbreak\vskip0pt plus.01\vsize\penalty-100
                \vskip0pt plus-.01\vsize
  	              \init{\titredeux}
                 \ifnum\optionparag=1{\init{\eqnumber}\init\enonce}\else{}\fi
                  \goodbreak
  	                {\parindent0mm\raggedcenter\fontetitreun\incr{\titreun}.\ 
                     \fontetitreun #1\par}\nobreak\medskip}
                     
\newtoks\titreunl
\titreunl={\ifnum\titreun=1{I}\fi%
\ifnum\titreun=2{II}\fi%
\ifnum\titreun=3{III}\fi%
\ifnum\titreun=4{IV}\fi%
\ifnum\titreun=5{V}\fi%
\ifnum\titreun=6{VI}\fi%
\ifnum\titreun=7{VII}\fi%
\ifnum\titreun=8{VIII}\fi%
\ifnum\titreun=9{IX}\fi%
\ifnum\titreun=10{X}\fi%
\ifnum\titreun=11{XI}\fi%
\ifnum\titreun=12{XII}\fi%
\ifnum\titreun=13{XIII}\fi%
}
\long\def\paraunl#1{\removelastskip\par\bigskip\bigskip\goodbreak\vskip0pt plus.01\vsize\penalty-100
                \vskip0pt plus-.01\vsize
  	              \init{\titredeux}\ifnum\optionparag=1{\init\eqnumber\init\enonce}\else{}\fi
                  \goodbreak{\fontetitreunl
  	                \Indentation{\global\advance\titreun by 1{\the\titreunl}.\ }{\fontetitreunl #1\par}}\nobreak\smallskip}

  
  \long\def\paradeux#1{\init{\titretrois}\vskip0pt plus.01\vsize\penalty-10
                \vskip0pt plus-.01\vsize\ifx \elie\oui\medskip\ifnum\titredeux>0\medskip\fi\fi
                 \Indentation{\fontetitredeux\the\titreun${\cdot}$\incr{\titredeux}.}
                              {\fontetitredeux\textfont1=\eleveni#1}\nobreak\par }
  
  \long\def\paradeuxb#1{\init{\titretrois}\vskip0pt plus.001\vsize\penalty-10
                \vskip0pt plus-.01\vsize{\ifx \elie\oui\medskip\ifnum\titredeux>0\medskip\fi\fi
                  \Indentation
  {\fontetitredeuxb\the\titreun${\cdot}$\incr{\titredeux}.}{ \fontetitredeuxb#1}}\nobreak
\smallskip}

\newtoks\titredeuxl
\titredeuxl={\ifnum\titredeux=1{A}\fi%
\ifnum\titredeux=2{B}\fi%
\ifnum\titredeux=3{C}\fi%
\ifnum\titredeux=4{D}\fi%
\ifnum\titredeux=5{E}\fi%
\ifnum\titredeux=6{F}\fi%
\ifnum\titredeux=7{G}\fi%
\ifnum\titredeux=8{H}\fi%
\ifnum\titredeux=9{I}\fi%
\ifnum\titredeux=10{J}\fi%
\ifnum\titredeux=11{K}\fi%
\ifnum\titredeux=12{L}\fi%
\ifnum\titredeux=13{M}\fi%
}
 \long\def\paradeuxl#1{\init{\titretrois}\vskip0pt plus.001\vsize\penalty-10
                \vskip0pt plus-.01
                \vsize \bigskip%
                  \Indentation
     {\fontetitredeuxl\global\advance\titredeux by 1
  \quad \the\titreunl${\cdot}$\the\titredeuxl.}{ \fontetitredeuxl#1}
  \removelastskip\nobreak\smallskip}
  

  \long\def\paratrois#1{\init{\titrequatre}\ifdim\lastskip<\smallskipamount
                \removelastskip\smallskip\fi
                 \vskip0pt plus.01\vsize\penalty-10
                  \vskip0pt
plus-.01\vsize{\ifx \elie\oui\ifnum\titretrois>0\medskip\fi\fi
\Indentation{\fontetitretrois\the\titreun${\cdot}$\the\titredeux${\cdot}$\incr{\titretrois}.\ }
  {\hskip0mm\baselineskip=14pt\fontetitretrois#1}\nobreak\smallskip}}
  
  
  \long\def\paratroisl#1{\init{\titrequatre}\ifdim\lastskip<\smallskipamount
                \removelastskip\fi
                 \vskip0pt plus.01\vsize\penalty-10
                  \vskip0pt
plus-.01\vsize\ifx \elie\oui\bigskip
\fi
\Indentation{\fontetitretrois\quad \quad \the\titreunl{${\cdot}$}\the\titredeuxl${\cdot}$\incr{\titretrois}.\ }
  {\hskip0mm\fontetitretrois#1}\nobreak\smallskip}


  \long\def\paraquatre#1{\ifdim\lastskip<\smallskipamount
                \removelastskip\smallskip\fi
                 \vskip0pt plus.01\vsize\penalty-10
                  \vskip0pt
                  plus-.01\vsize\par
 
\Indentation{\fontetitrequatre \the\titreun{${\cdot}$}\the\titredeux${\cdot}$\the\titretrois${\cdot}$\incr{\titrequatre}.\ }
{\hskip0mm\fontetitrequatre#1}\nobreak\smallskip}


\newtoks\titrequatrel
\titrequatrel={\ifnum\titrequatre=1{a}\fi%
\ifnum\titrequatre=2{b}\fi%
\ifnum\titrequatre=3{c}\fi%
\ifnum\titrequatre=4{d}\fi%
\ifnum\titrequatre=5{e}\fi%
\ifnum\titrequatre=6{f}\fi%
\ifnum\titrequatre=7{g}\fi%
\ifnum\titrequatre=8{h}\fi%
\ifnum\titrequatre=9{i}\fi%
\ifnum\titrequatre=10{j}\fi%
\ifnum\titrequatre=11{k}\fi%
\ifnum\titrequatre=12{l}\fi%
\ifnum\titrequatre=13{m}\fi%
}
\long\def\paraquatrel#1{\ifdim\lastskip<\smallskipamount
                \removelastskip\smallskip\fi
                 \vskip0pt plus.01\vsize\penalty-10
                  \vskip0pt
                  plus-.01\vsize{\bigskip
\Indentation{\global\advance\titrequatre by 1
\fontetitrequatre\quad \quad \quad \the\titreunl${\cdot}$\the\titredeuxl${\cdot}$\the\titretrois${\cdot}$\the\titrequatrel.\ }
{\hskip0mm\fontetitrequatre#1}\nobreak\smallskip}}

\ifx\optionkeys\oui
\def\drefun#1{\definexref{¤#1}{{\the\titreun}}{}} 
\def\drefdeux#1{\definexref{¤#1}{{\the\titreun}.{\the\titredeux}}{}}
\def\dreftrois#1{\definexref{¤#1}{{\the\titreun}.{\the\titredeux}.{\the\titretrois}}{}}
\else
\def\drefun#1{\definexref{prg#1}{{\the\titreun}}{}} 
\def\drefdeux#1{\definexref{prg#1}{{\the\titreun}.{\the\titredeux}}{}}
\def\dreftrois#1{\definexref{prg#1}{{\the\titreun}.{\the\titredeux}.{\the\titretrois}}{}}
\fi

%


  \long\def\propdeux#1#2#3#4{%
       \avance{\enonce}
       \leavevmode\edef\temp{#2}%
         \ifx\temp\empty 
          \else
           \definexref{#2}{#1~{\the\titreun.\the\enonce}}{enonces}
            \definexref{s#2}{{\the\titreun.\the\enonce}}{enonces}
             \fi
\smallskip
      \noindent{\bf#1\ {\bf\the\titreun.\the\enonce{#3}.}\enspace}{\sl#4\par}%
      \ifdim\lastskip<\medskipamount \removelastskip\penalty55\par \fi
   }

  \long\def\propun#1#2#3#4{%
      \avance{\enonce}
       \leavevmode\edef\temp{#2}%
        \ifx\temp\empty 
          \else
           \definexref{#2}{#1~{\the\enonce}}{enonces}
            \definexref{{s#2}}{{\the\enonce}}{enonces}
             \fi
   \par 
     \noindent{\bf#1\ {\bf\the\enonce{#3}.}\enspace}{\sl#4\par}%
     \ifdim\lastskip<\medskipamount \removelastskip\penalty55\medskip\fi
  }
  
  \long\def\prop#1#2#3#4{\ifnum\optionparag=1
                          \propdeux{#1}{#2}{\textfont1=\elevenib#3}{#4} \else\propun{#1}{#2}{\textfont1=\elevenib#3}{#4}\fi}

  \long\def\propt#1#2#3{\ifx\tpf\oui \prop{Th\'eo\-r\`eme}{#1}{#2}{#3}\par
                       \else\prop{Theorem}{#1}{#2}{#3}\par\fi}
  \long\def\Propt#1#2{\propt{#1}{}{#2}}
  \long\def\propl#1#2#3{\ifx\tpf\oui\prop{Lem\-me}{#1}{#2}{#3}\par
                         \else\prop{Lemma}{#1}{#2}{#3}\par\fi}
  \long\def\Propl#1#2{\propl{#1}{}{#2}}
  \long\def\propc#1#2#3{\ifx\tpf\oui\prop{Corol\-laire}{#1}{#2}{#3}\par
                         \else\prop{Corollary}{#1}{#2}{#3}\par\fi}
  \long\def\Propc#1#2{\propc{#1}{}{#2}}

  \long\def\propd#1#2#3{\ifx\tpf\oui\prop{D\'efi\-nition}{#1}{#2}{#3}\par
                       \else\prop{Definition}{#1}{#2}{#3}\par\fi} 
  
  \long\def\proptd#1#2#3{\ifx\tpf\oui\prop{Th\'eor\`eme et d\'efi\-nition}{#1}{#2}{#3}\par
                       \else\prop{Theorem and definition}{#1}{#2}{#3}\par\fi}


  
  \newcount\optionparag\optionparag=1
  
  \long\def\section#1#2{\ifnum\optionparag=1 \paraun{#2} 
                        \else\goodbreak{\fontetitreun
  	                \Indentation{#1.\ }{#2}}\nobreak\smallskip\fi}

  \def\eqconstruct#1{\ifnum\optionparag=1{\the\titreun\hbox{$\cdot$}#1}\else{#1}\fi}

  
  
  \def\numref{oui}  
  
  \newcount\mesref\mesref=0 
  \def\defbib#1{\ifx\numref\oui\global\advance\mesref by 1 \definexref{#1}{{\the
                 \mesref}}{}\else\definexref{#1}{#1}{}\fi}
  \def\bibtem#1{\defbib{#1}\item{\citer{#1}}}
  \def\citer#1{[\ref{#1}]}

  
  \font\seventeenmsa=msam10 at 17pt    
  \font\fourteenmsa=msam10 at 14pt
  \font\twelvemsa=msam10 at 12pt
  \font\tenmsa=msam10                 
  \font\ninemsa=msam10 at 9pt 
  \font\eightmsa=msam10 at 8pt 
  \font\sevenmsa=msam7 
  \font\sixmsa=msam10 at 6pt
  \font\fivemsa=msam5
  \newfam\msafam\textfont\msafam=\tenmsa\scriptfont\msafam=\sevenmsa\scriptscriptfont\msafam=\fivemsa
  
  \font\seventeenbb=msbm10 at 17pt     
  \font\fourteenbb=msbm10 at 14pt
  \font\twelvebb=msbm10 at 12pt
  \font\tenbb=msbm10                   
  \font\ninebb=msbm10 at 9pt 
  \font\eightbb=msbm10 at 8pt 
  \font\sevenbb=msbm7 
  \font\sixbb=msbm10 at 6pt
  \font\fivebb=msbm5 
  \newfam\bbfam\textfont\bbfam=\tenbb\scriptfont\bbfam=\sevenbb\scriptscriptfont\bbfam=\fivebb
  \def\bb{\fam\bbfam\tenbb}%

  \font\seventeenscaln=eusm10 at 17pt   
  \font\twelvescaln=eusm10 at 12pt
  \font\tenscaln=eusm10                
  \font\ninescaln=eusm10 scaled 900
  \font\eightscaln=eusm10 scaled 800
  \font\sevenscaln=eusm10 scaled 700
  \font\sixscaln=eusm10 scaled 600
   
  \newfam\scalnfam\textfont\scalnfam=\tenscaln\scriptfont\scalnfam=\sevenscaln\scriptscriptfont\scalnfam=\sixscaln
  \def\scaln{\fam\scalnfam\tenscaln}%
  \def\scal{\scaln}
  
  \font\tenscalb=eusb10                

  \font\sevenscalb=eusb10 scaled 700

  \newfam\scalbfam\textfont\scalbfam=\tenscalb\scriptfont\scalbfam=\sevenscalb
  %
  
  %
  %
  \font\fourteenrm=cmr12 scaled 1200
  \font\elevenrm=cmr10 at 11pt
  \font\twelverm=cmr12
  \font\ninerm=cmr9
  \font\eightrm=cmr8      
  \font\sevenrm=cmr7
  \font\sixrm=cmr6

  \font\seventeenpcap=cmcsc10 at 17pt
  \font\tenpcap=cmcsc10                        
  \font\ninepcap=cmcsc9
  \font\eightpcap=cmcsc8
  \font\sevenpcap=cmcsc10 scaled 700
  
  \newfam\pcapfam\textfont\pcapfam=\tenpcap\scriptfont\pcapfam=\sevenpcap
  \def\pcap{\fam\pcapfam\tenpcap}
  
  \font\seventeenrm=cmbx12 scaled 1400

  \font\fourteenbf=cmbx10 scaled 1400
  
  \font\twelvebf=cmbx12
  \font\elevenbf=cmbx10 at 11pt
  \font\ninebf=cmbx9  
  \font\eightbf=cmbx8
  \font\sixbf=cmbx6
  
  \font\tengot=eufm10                           
   
  \font\eightgot=eufm10 at 8truept 
  \font\sevengot=eufm7 
  \font\sixgot=eufm10 at 6 truept 
   
  \newfam\gotfam
  \textfont\gotfam=\tengot\scriptfont\gotfam=\sevengot\scriptscriptfont\gotfam=\sixgot
  \def\got{\fam\gotfam\tengot}%

  
  \def\tit{%
  \textfont0=\seventeenrm\scriptfont0=\tenrm\def\rm{\fam0\seventeenrm}%
  \textfont1=\seventeenib\scriptfont1=\twelveib%
  \textfont2=\seventeensy\scriptfont2=\twelvesy\scriptscriptfont2=\ninesy
  \textfont3=\seventeenex
  \textfont\itfam=\seventeenti
  \def\it{\fam\itfam\seventeenti}%
  \textfont\bbfam=\seventeenbb \scriptfont\bbfam=\twelvebb
  \def\bb{\fam\bbfam\seventeenbb}%
  \textfont\msafam=\seventeenmsa\scriptfont\msafam=\twelvemsa
  \textfont\scalnfam=\seventeenscaln
  \def\pcap{\fam\pcapfam\seventeenpcap}
  \normalbaselineskip=25pt\normalbaselines\rm}

  \font\seventeenti=cmbxti10 scaled 1680
  
  \font\fourteenti=cmbxti10 at 14pt
  
  \font\twelveti=cmbxti10 scaled 1200
  \font\eleventi=cmbxti10 at 11pt

  %
  %
  \font\twelveit=cmti12	
  \font\elevenit=cmti10 scaled 1100
  \font\nineit=cmti9
  \font\eightit=cmti8
  \font\sevenit=cmti7

  %
  %
 
 \font\seventeenib=cmmib10 scaled 1680
  \font\fourteenib=cmmib10 scaled 1400
  \font\twelveib=cmmib10 scaled 1200
  \font\elevenib=cmmib10 scaled 1100
  \font\tenib=cmmib10
\font\eightib=cmmib10 scaled 800
  \font\nineib=cmmib10 scaled 900
\font\sevenib=cmmib10 scaled 700
\font\sixib=cmmib10 scaled 600
\font\fiveib=cmmib10 scaled 500

\ifx\ITAN\oui
\else
\innernewfam\cmmibfam
\textfont\cmmibfam=\tenib
\scriptfont\cmmibfam=\sevenib
\scriptscriptfont\cmmibfam=\fiveib
\def\ib{\fam\cmmibfam\tenib}
\fi

  %
  %
  
  \font\eleveni=cmmi10 scaled 1100
  \font\ninei=cmmi9
  \font\eighti=cmmi8 
  \font\seveni=cmmi7 			                
  \font\sixi=cmmi6
  
  \font\ninesl=cmsl9                    
  \font\eightsl=cmsl8 
  \font\sevensl=cmsl10 at 7pt

  \font\ninett=cmtt9                    
  \font\eighttt=cmtt8
  \font\seventt=cmtt10 scaled 700

  \font\seventeensy=cmsy10 scaled 1680    
  \font\fourteensy=cmsy10 scaled 1400
  \font\twelvesy=cmsy10 scaled 1176
  
  \font\ninesy=cmsy9                      
  \font\eightsy=cmsy8
  \font\sixsy=cmsy6
  \font\seventeenex=cmex10 at 17pt
  \font\fourteenex=cmex10 at 14pt
  \font\twelveex=cmex10 at 12pt
  \font\nineex=cmex10 at 9pt
  \font\eightex=cmex10 at 8pt
  \font\sevenex=cmex10 at 7pt
  \font\sixex=cmex10 at 6pt
  \font\fiveex=cmex10 at 5pt
  
   
  \font\fourteengp=cmmi10 at 14pt
  
  \font\twelvegp=cmmib10 at 12pt
  \font\elevengp=cmmib10 at 11pt
  \font\tengp=cmmib10                          
  \font\ninegp=cmmib10 at 9pt 
  \font\eightgp=cmmib8 
   
  \font\sixgp=cmmib6


  \def\gponze{\textfont0=\elevenbf\scriptfont0=\eightbf\scriptscriptfont0=\sixbf
           \textfont1=\elevengp\scriptfont1=\eightgp\scriptscriptfont1=\sixgp}
  \def\gpdouze{\textfont0=\twelvebf\scriptfont0=\tenbf\scriptscriptfont0=\ninebf
           \textfont1=\twelvegp\scriptfont1=\tengp\scriptscriptfont1=\ninegp}        
  
 \def\gpquatorze{\textfont0=\fourteenbf\scriptfont0=\twelvebf\scriptscriptfont0=\elevenbf
           \textfont1=\fourteengp\scriptfont1=\twelvegp\scriptscriptfont1=\elevengp}

  
  \expandafter\chardef\csname pre amssym.def at\endcsname=\the\catcode`\@
  \catcode`\@=11
  \def\undefine#1{\let#1\undefined}
  \def\newsymbol#1#2#3#4#5{\let\next@\relax
   \ifnum#2=\@ne\let\next@\msafam@\else
   \ifnum#2=\tw@\let\next@\bbfam@\fi\fi
   \mathchardef#1="#3\next@#4#5}
  \def\mathhexbox@#1#2#3{\relax
   \ifmmode\mathpalette{}{\m@th\mathchar"#1#2#3}%
   \else\leavevmode\hbox{$\m@th\mathchar"#1#2#3$}\fi}
  \def\hexnumber@#1{\ifcase#1 0\or 1\or 2\or 3\or 4\or 5\or 6\or 7\or 8\or
   9\or A\or B\or C\or D\or E\or F\fi}
  
  \def\setboxz@h{\setbox\z@\hbox}
  \def\wdz@{\wd\z@}
  \def\boxz@{\box\z@}
  
  \edef\msafam@{\hexnumber@\msafam}
  \mathchardef\dabar@"0\msafam@39
  
  \edef\bbfam@{\hexnumber@\bbfam}
  \def\widehat#1{\setboxz@h{$\m@th#1$}%
   \ifdim\wdz@>\tw@ em\mathaccent"0\bbfam@5B{#1}%
   \else\mathaccent"0362{#1}\fi}
  \def\widetilde#1{\setboxz@h{$\m@th#1$}%
   \ifdim\wdz@>\tw@ em\mathaccent"0\bbfam@5D{#1}%
   \else\mathaccent"0365{#1}\fi}
  \newsymbol\leqq 1335          
  \newsymbol\leqslant 1336
  \newsymbol\lessgtr 1337       
  \newsymbol\backprime 1038     
  \newsymbol\risingdotseq 133A  
  \newsymbol\fallingdotseq 133B 
  \newsymbol\succcurlyeq 133C   
  \newsymbol\geqq 133D          
  \newsymbol\geqslant 133E
  \newsymbol\nmid 232D
  \newsymbol\nexists 2040
  \newsymbol\smallsetminus 2272
  \newsymbol\varnothing 203F
  
  \catcode`\@=\active

  \catcode`\@=11
  \newcount\typofr\typofr=1
  
  \catcode`\;=\active
  \def;{\ifnum\typofr=1\relax\ifhmode\ifdim\lastskip>\z@\unskip\fi
     \kern.2em\fi\string;\else\string;\fi}
  
  \catcode`\:=\active
  \def:{\ifnum\typofr=1\relax\ifhmode\ifdim\lastskip>\z@\unskip\fi
  \penalty\@M\ \fi\string:\else\string:\fi}
  
  \catcode`\!=\active
  \def!{\ifnum\typofr=1\relax\ifhmode\ifdim\lastskip>\z@\unskip\fi
     \kern.2em\fi\string!\else\string!\fi}
  
  \catcode`\?=\active
  \def?{\ifnum\typofr=1\relax\ifhmode\ifdim\lastskip>\z@\unskip\fi
     \kern.2em\fi\string?\else\string?\fi}

  \def\francais{\typofr=1\def\tpf{oui}}
  
  \def\oui{oui}
  \francais
  
  \catcode`\@=12
  

\ifx\textures\oui
\def\raye #1|{\leavevmode\setbox1=\hbox{#1}%
\raise .5pt\hbox to \wd1{\xleaders\hbox{\rge{$ \sct / $}%
\kern 1pt}\hfill\kern -1pt }\kern -\wd1 \unhbox1\relax }

\def\barre #1|{\leavevmode\setbox1=\hbox{#1}%
\rlap{\Red\vrule height 2.4pt depth -1.2pt width \wd1}\Black \unhbox1\relax}
\else
\def\raye #1|{\leavevmode\setbox1=\hbox{#1}%
\raise .5pt\hbox to \wd1{\xleaders\hbox{\rge{$ \sct / $}%
\kern 1pt}\hfill\kern -1pt }\kern -\wd1 \unhbox1\relax }

\def\barre #1|{\leavevmode\setbox1=\hbox{#1}%
\rlap{\color{red}\vrule height 2.4pt depth -1.2pt width \wd1}\color{black} \unhbox1\relax}

\fi
  

  
  \def\og{\leavevmode\raise.24ex\hbox{$\scriptscriptstyle\langle\!\langle\>$}}    
  \def\fg{\leavevmode\raise.24ex\hbox{$\scriptscriptstyle\>\rangle\!\rangle$}}    

  \def\d{\,{\rm d}}
  \def\dd{{\rm d}}

  \def\r{{\bb R}}
  \def\CC{{\bb C}}
  \def\N{{\bb N}}

  \def\PP{{\bb P}}

  \def\A{{\scal A}}

  \def\HH{{\scal H}}

  \def\L{{\scal L}}
  \def\M{{\scal M}}
  
  \def\O{{\scal O}}
  \def\P{{\scaln P}}

  \def\T{{\scal T}}

  \def\frac#1#2{{#1\over #2}}
  \def\di#1#2{\sct#1\atop{\sct#2}}
  \def\tri#1#2#3{{\sct#1\atop\sct#2}\atop\sct#3}

  \def\numero{n$^{\rm o}\thinspace$}
\def\numeros{n$^{\rm os}\thinspace$}

  \def\qedbox{$\rlap{$\sqcap$}\sqcup$}           
  \def\qed{\nobreak\hfill\penalty250 \hbox{}\nobreak\hfill\qedbox\par }

  \def\sumast{\mathop{{\sum}^*}}
  
  \def\numero{n$^{\rm o}\thinspace$}
  \def\np{nombre premier}
  \def\nps{nombres premiers}

  \def\¤{\S\thinspace}

  \def\¥{$\bullet$ }
  
  
  \def\e{{\rm e}}
  \def\mod{\mathop{\rm mod}\nolimits}
  \def\md#1#2{\equiv#1\,({\rm mod\,}#2)}
  \def\no#1{\Vert#1\Vert}

  \def\epsilon{\varepsilon}

  \def\phi{\varphi}
  \def\theta{\vartheta}
  \def\rho{\varrho}
  \def\dm{{\textstyle{1\over 2}}}
  \def\txt{\textstyle}
  \def\dsp{\displaystyle}
  \def\sct{\scriptstyle}
  \def\pf{\noi{\it Proof. }}
  \def\nid{\ifnum\typofr=1\par\noindent{\it D\'emonstration. }\else\pf\fi}
  \def\noi{\noindent}
  \def\rem{\ifnum\typofr=1\noi{\it Remarque.}\ \else\noi{\it Remark.}\ \fi}
  \def\rems{\ifnum\typofr=1\noi{\it Remarques.}\ \else\noi{\it Remarks.}\ \fi}
  \def\re{{\Re e\,}}
  \def\im{{\Im m\,}}
  \def\ov{\overline}

  \def\emptyset{{\varnothing}}

  \def\1{{\bf 1}}
  \def\|{\Vert}

  \def\leq{\leqslant}
  \def\geq{\geqslant}
  
  \def\cf{{cf.}}


  \def\fl#1{\left\lfloor #1 \right\rfloor}

  \def\log{\mathop{\rm log}\nolimits}
  \def\ft#1#2{{\txt{#1\over #2}}}


  
\def\vbs#1{\left|#1\right|}

\def\abs#1{\left|#1\right|}


  \def\pmb#1{\setbox0=\hbox{#1}%
  \kern-.025em\copy0\kern-\wd0\kern.05em\copy0\kern-\wd0\kern-.025em\raise .0433em\box0 }

  
  \skewchar\eighti='177 \skewchar\sixi='177
  \skewchar\eightsy='60 \skewchar\sixsy='60
  
  \def\eightpoint{%
  \textfont0=\eightrm\scriptfont0=\sixrm\scriptscriptfont0=\fiverm
  \def\rm{\fam0\eightrm}%
  \textfont1=\eighti\scriptfont1=\sixi
  \scriptscriptfont1=\fivei\def\oldstyle{\fam1\seveni}%
  \textfont2=\eightsy\scriptfont2=\sixsy\scriptscriptfont2=\fivesy
  \textfont3=\eightex\scriptfont3=\sixex
  \textfont\itfam=\eightit
  \def\it{\fam\itfam\eightit}%
  \textfont\slfam=\eightsl
  \def\sl{\fam\slfam\eightsl}%
  \textfont\bbfam=\eightbb \scriptfont\bbfam=\sixbb\scriptscriptfont\bbfam=\fivebb
  \def\bb{\fam\bbfam\eightbb}%
  \textfont\msafam=\eightmsa\scriptfont\msafam=\sixmsa
  \textfont\scalnfam=\eightscaln
  \def\scaln{\fam\scalnfam\eightscaln}
  \textfont\ttfam=\eighttt
  \def\tt{\fam\ttfam\eighttt}%
\textfont\gotfam=\eightgot
  \textfont\bffam=\eightbf\scriptfont\bffam=\sixbf\scriptscriptfont\bffam=\fivebf
  \def\bf{\fam\bffam\eightbf}%
  \ifx\ITAN\oui\else\textfont\cmmibfam=\eightib
       \scriptfont\cmmibfam=\sixib
        \scriptscriptfont\cmmibfam=\fiveib
         \def\ib{\fam\cmmibfam\eightib}
   \fi
  \textfont\pcapfam=\eightpcap
  \def\pcap{\fam\pcapfam\eightpcap}
  \abovedisplayskip=2pt plus2pt minus 2pt
  \belowdisplayskip=2pt plus1pt minus 2pt
  \abovedisplayshortskip= 1pt plus 2pt minus 1pt
  \belowdisplayshortskip= 1pt plus 2pt minus 1pt
  \smallskipamount=2pt plus 1pt minus 2pt
  \medskipamount=3pt plus 2pt minus 2pt
  \bigskipamount=7pt plus 3pt minus 3pt
  \setbox\strutbox=\hbox{\vrule height 5pt depth 2pt width 0pt}%
  \normalbaselineskip=9pt\normalbaselines\rm}

  \def\({\left(}
  \def\){\right)}
  
  \def\footnoterule{\kern -2pt\hrule width 7truecm\kern 2.4pt}
  
  \def\xnotedef#1{\definexref{#1}{\noexpand\number\footnotenumber}{Note}}%

  
  
  \def\ninepoint{%
  \textfont0=\ninerm\scriptfont0=\sixrm\scriptscriptfont0=\fiverm
  \def\rm{\fam0\ninerm}%
  \textfont1=\ninei\scriptfont1=\sixi
  \scriptscriptfont1=\fivei\def\oldstyle{\fam1\ninei}%
  \textfont2=\ninesy\scriptfont2=\sixsy\scriptscriptfont2=\fivesy
  \textfont3=\nineex\scriptfont3=\sixex
  \textfont\itfam=\nineit
  \def\it{\fam\itfam\nineit}%
  \textfont\slfam=\ninesl
  \def\sl{\fam\slfam\ninesl}%
  \textfont\bbfam=\ninebb\scriptfont\bbfam=\sixbb\scriptscriptfont\bbfam=\fivebb
  \def\bb{\fam\bbfam\ninebb}%
  \textfont\msafam=\ninemsa\scriptfont\msafam=\sixmsa\scriptscriptfont\msafam=\fivemsa
  \textfont\scalnfam=\ninescaln
  \def\scaln{\fam\scalnfam\ninescaln}
  \textfont\ttfam=\ninett
  \def\tt{\fam\ttfam\ninett}%
  \textfont\bffam=\ninebf\scriptfont\bffam=\sixbf\scriptscriptfont\bffam=\fivebf
  \def\bf{\fam\bffam\ninebf}%
  \abovedisplayskip=3pt plus2pt minus 2pt
  \belowdisplayskip=3pt plus1pt minus 2pt
  \abovedisplayshortskip= 2pt plus 2pt minus 1pt
  \belowdisplayshortskip= 2pt plus 2pt minus 1pt
  \smallskipamount=2pt plus 1pt minus 2pt
  \medskipamount=3pt plus 2pt minus 2pt
  \bigskipamount=7pt plus 3pt minus 3pt
  \setbox\strutbox=\hbox{\vrule height 5pt depth 2pt width 0pt}%
  \normalbaselineskip=10.5pt plus.3pt minus.3pt\normalbaselines\rm}

  \def\sevenpoint{%
  \textfont0=\sevenrm\scriptfont0=\sixrm\scriptscriptfont0=\fiverm
  \def\rm{\fam0\sevenrm}%
  \textfont1=\seveni\scriptfont1=\sixi
  \scriptscriptfont1=\fivei\def\oldstyle{\fam1\seveni}%
  \textfont2=\sevensy\scriptfont2=\sixsy\scriptscriptfont2=\fivesy
  \textfont3=\sevenex\scriptfont3=\fiveex
  \textfont\itfam=\sevenit
  \def\it{\fam\itfam\sevenit}%
  \textfont\slfam=\sevensl
  \def\sl{\fam\slfam\sevensl}%
  \textfont\bbfam=\sevenbb \scriptfont\bbfam=\sixbb\scriptscriptfont\bbfam=\fivebb
  \def\bb{\fam\bbfam\sevenbb}%
  \textfont\msafam=\sevenmsa\scriptfont\msafam=\sixmsa
  \textfont\scalnfam=\sevenscaln
  \def\scaln{\fam\scalnfam\sevenscaln}
  \textfont\bffam=\sevenbf\scriptfont\bffam=\sixbf\scriptscriptfont\bffam=\fivebf
  \def\bf{\fam\bffam\sevenbf}%
  \textfont\ttfam=\seventt
  \abovedisplayskip=2pt plus2pt minus 2pt
  \belowdisplayskip=2pt plus1pt minus 2pt
  \abovedisplayshortskip= 1pt plus 2pt minus 1pt
  \belowdisplayshortskip= 1pt plus 2pt minus 1pt
  \smallskipamount=2pt plus 1pt minus 2pt
  \medskipamount=3pt plus 2pt minus 2pt
  \bigskipamount=7pt plus 3pt minus 3pt
  \setbox\strutbox=\hbox{\vrule height 5pt depth 2pt width 0pt}%
  \normalbaselineskip=9pt\normalbaselines\rm}

 \def\onzepts{%
 \textfont0=\elevenrm\scriptfont0=\ninerm
 \textfont1=\elevenib\scriptfont1=\ninei}

\def\douzepts{%
  \textfont0=\twelverm\scriptfont0=\tenrm\def\rm{\fam0\twelverm}%
  \textfont1=\twelveib\scriptfont1=\teni%
  \textfont2=\twelvesy\scriptfont2=\tensy\scriptscriptfont2=\eightsy
  \textfont3=\twelveex
  \textfont\itfam=\twelveti
  \def\it{\fam\itfam\twelveti}%
  \textfont\bffam=\twelvebf\scriptfont\bffam=\tenbf\scriptscriptfont\bffam=\eightbf
  \def\bf{\fam\bffam\twelvebf}%
  \textfont\bbfam=\twelvebb \scriptfont\bbfam=\tenbb
  \def\bb{\fam\bbfam\twelvebb}%
  \textfont\msafam=\twelvemsa\scriptfont\msafam=\tenmsa
  \textfont\scalnfam=\twelvescaln
  \normalbaselineskip=15pt\normalbaselines\rm}

\def\quatorzepts{%
  \textfont0=\fourteenrm\scriptfont0=\twelverm\def\rm{\fam0\fourteenrm}%
  \textfont1=\fourteenib\scriptfont1=\twelveib%
  \textfont2=\fourteensy\scriptfont2=\twelvesy\scriptscriptfont2=\tensy
  \textfont3=\fourteenex
  \textfont\itfam=\fourteenti
  \def\it{\fam\itfam\fourteenti}%
  \textfont\bffam=\fourteenbf\scriptfont\bffam=\twelvebf\scriptscriptfont\bffam=\tenbf
  \def\bf{\fam\bffam\fourteenbf}%
  \textfont\bbfam=\fourteenbb \scriptfont\bbfam=\twelvebb
  \def\bb{\fam\bbfam\fourteenbb}%
  \textfont\msafam=\fourteenmsa\scriptfont\msafam=\twelvemsa
  \textfont\scalnfam=\twelvescaln
  \normalbaselineskip=18pt\normalbaselines\rm}


\def\AA{{\it Acta Arith.}}

\def\picture #1 by #2 (#3){\leavevmode\vbox to #2{
     \hrule width #1 height 0pt depth 0pt
      \vfill
       \special{picture #3}}}

\def\illustration #1 by #2 (#3) scaled #4{\dimen1=#2
  \divide\dimen1 by 1000
  \multiply\dimen1 by #4
  \vtop to \dimen1{\dimen1=#1
  \divide\dimen1 by 1000
  \multiply\dimen1 by #4
  \hsize=\dimen1\vss
  \noindent\special{illustration #3 scaled #4}}}

 \fi

\dimstand  
 \francais
\vsize=250truemm
\voffset-5mm

\ifx\optionkeymacros\undefined\else \fi

\catcode`\Œ=\active\defŒ{{\aa}}       
\catcode`\º=\active\defº{\int}        
\catcode`\=\active\def{\c c}        
\catcode`\¶=\active\def¶{\partial}    
\catcode`\Ä=\active\defÄ{\oint}       
\catcode`\Æ=\active\defÆ{\triangle}   
\catcode`\Â=\active\defÂ{\neg}        
\catcode`\µ=\active\defµ{\mu}         
\catcode`\¿=\active\def¿{{\o}}        
\catcode`\¹=\active\def¹{\pi}         
\catcode`\Ï=\active\defÏ{{\oe}}       
\catcode`\§=\active\def§{{\ss}}       
\catcode`\ =\active\def {\dagger}     
\catcode`\Ã=\active\defÃ{\sqrt}       
\catcode`\·=\active\def·{\Sigma}      
\catcode`\Å=\active\defÅ{\approx}     
\catcode`\½=\active\def½{\Omega}      
\catcode`\£=\active\def£{{\it\$}}     
\catcode`\°=\active\def°{\infty}      
\catcode`\¤=\active\def¤{{\S}}        
\catcode`\¦=\active\def¦{{\P}}        
\catcode`\¥=\active\def¥{\bullet}     
\catcode`\»=\active\def»{\leavevmode\raise.585ex\hbox{\b a}}      
\catcode`\¼=\active\def¼{\leavevmode\raise.6ex\hbox{\b o}}        
\catcode`\­=\active\def­{\not=}       
\catcode`\²=\active\def²{\leq}        
\catcode`\³=\active\def³{\geq}        
\catcode`\Ö=\active\defÖ{\div}        
\catcode`\É=\active\defÉ{{\dots}}     
\catcode`\¾=\active\def¾{{\ae}}       
\catcode`\Ç=\active\defÇ{\og}         
\catcode`\Ò=\active\defÒ{``}          
\catcode`\Á=\active\defÁ{!`}          
\catcode`\¢=\active\def¢{\rlap/c}     
\catcode`\Ô=\active\defÔ{`}           
\catcode`\Õ=\active\defÕ{'}           


\catcode`\=\active\def{{\AA}}       
\catcode`\'=\active\def'{\c C}        
\catcode`\¯=\active\def¯{{\O}}        
\catcode`\¸=\active\def¸{\Pi}         
\catcode`\Î=\active\defÎ{{\OE}}       
\catcode`\®=\active\def®{{\AE}}       
\catcode`\×=\active\def×{\diamond}    
\catcode`\¡=\active\def¡{\accent'27}  
\catcode`\Ó=\active\defÓ{''}          
\catcode`\±=\active\def±{\pm}         
\catcode`\È=\active\defÈ{\fg}         
\catcode`\À=\active\defÀ{?`}          
\catcode`\Ð=\active\defÐ{--}          
\catcode`\Ñ=\active\defÑ{---}         


\catcode`\Š=\active\defŠ{\"a}        
\catcode`\'=\active\def'{\"e}        
\catcode`\•=\active\def•{\"{\i}}     
\catcode`\š=\active\defš{\"o}        
\catcode`\Ÿ=\active\defŸ{\"u}        
\catcode`\Ø=\active\defØ{\"y}        
\catcode`\å=\active\defå{\^A}        
\catcode`\€=\active\def€{\"A}        
\catcode`\…=\active\def…{\"O}        
\catcode`\†=\active\def†{\"U}        
\catcode`\‡=\active\def‡{\'a}        
\catcode`\Ž=\active\defŽ{\'e}        
\catcode`\'=\active\def'{\'{\i}}     
\catcode`\—=\active\def—{\'o}        
\catcode`\œ=\active\defœ{\'u}        
\catcode`\ƒ=\active\defƒ{\'E}        
\catcode`\æ=\active\defæ{\^E}        
\catcode`\é=\active\defé{\`E}        %
\catcode`\ˆ=\active\defˆ{\`a}        
\catcode`\=\active\def{\`e}        
\catcode`\"=\active\def"{\`{\i}}     
\catcode`\˜=\active\def˜{\`o}        
\catcode`\=\active\def{\`u}        
\catcode`\Ë=\active\defË{\`A}        
\catcode`\‹=\active\def‹{\~a}        
\catcode`\–=\active\def–{\~n}        
\catcode`\›=\active\def›{\~o}        
\catcode`\Ì=\active\defÌ{\~A}        
\catcode`\"=\active\def"{\~N}        
\catcode`\Í=\active\defÍ{\~O}        
\catcode`\‰=\active\def‰{\^a}        
\catcode`\=\active\def{\^e}        
\catcode`\"=\active\def"{\^{\i}}     
\catcode`\™=\active\def™{\^o}        
\catcode`\ž=\active\defž{\^u}        

\let\optionkeymacros\null

\def\paradouze{oui}


\def\voircorrections{oui} 

\ifx\montrerlabels\oui\input montrerlabels.tex\fi

%
%

\newcount\paras\paras=0
\newcount\sparas\sparas=0
\def\tocsectionentry#1#2{\init{\sparas}\avance\paras
                          {\quad\bf\the\paras\quad }{\hskip-2mm
#1\dotfill\hskip3mm\rm#2}\par}%
\def\tocsubsectionentry#1#2{\avance\sparas
                          {\qquad\eightpoint\it\the\paras.\the\sparas}
{\hskip-2mm\eightpoint\it#1\dotfill\hskip3mm\rm#2}\par}%

%
%

\def\??{\leftnote{{??}}}

\ifx\voircorrections\oui\else\def\rge#1{#1}\fi 

\dateheure


\font\tenib=cmmib10
            \font\nineib=cmmib10 scaled 900
\font\sevenib=cmmib10 scaled 700
\font\eightib=cmmib10 scaled 800
\font\fiveib=cmmib10 scaled 500
\def\itg{\ib}

\def\ITAN{oui}

    \font\tenrsfs=rsfs7 at 10pt

    \font\sevenrsfs=rsfs7
    \font\sixrsfs=rsfs7 at 6pt
    
\newfam\rsfsfam\textfont\rsfsfam=\tenrsfs\scriptfont\rsfsfam=\sevenrsfs\scriptscriptfont\rsfsfam=\sixrsfs
    \def\rsfs{\fam\rsfsfam\tenrsfs}%
\def\L{{\rsfs L}}
\def\T{{\scal T}}

\def\bfc{{\itg c}} 
\def\bfr{{\itg r}}

\def\gotS{{\got S}}
\def\gotQ{{\got Q}}
\def\gotZ{{\got Z}}

\def\fl#1{\left\lfloor #1 \right\rfloor}

\def\titrart{Sommes de G‡l et applications\anote{*}{Nous incluons ici certaines corrections mineures relativement ˆ la version publiŽe.\hfill}}
\def\auteurs{RŽgis de la Bret\`eche \& GŽrald Tenenbaum}
\hautspages{\auteurs}{\titrart}
\titrecentre{\titrart}
\bigskip
\centerline{\auteurs}
\bigskip\bigskip
{\eightpoint\leftskip1cm\rightskip1cm
\noi{\bf Abstract.} 
We evaluate the asymptotic size of various sums of G‡l type, in particular
$$S( \M):=\sum_{m,n\in\M} \sqrt{(m,n) \over [m,n]},$$
where $\M$ is a finite set of integers.
Elaborating on methods  recently developed by Bondarenko and Seip, we obtain  an asymptotic formula   for $$\log\Big(
\sup_{|\M|= N}{S( \M)/N}\Big)$$
and derive new lower bounds for localized extreme values of the Riemann zeta-function, for extremal values  of  some Dirichlet $L$-functions at $s=\dm$, and for large character sums. 
\PAR
\medskip\noi
{ \bf Keywords:}  G‡l's theorem, GCD sums, Dirichlet polynomials, the Riemann zeta function, Dirichlet $L$-functions, character sums, resonance method.\par 
\medskip\noi
{\bf AMS 2010 Classification:} 11A05,  11L40, 11M06, 11N25, 11N37.\par}
\par 
\bigskip\bigskip

\medskip 
\paraun{Introduction et ŽnoncŽ des rŽsultats} 
 \paradeux{Sommes de G‡l} 
Soit $\M$ un ensemble fini de    $N$ entiers.  Nous considŽrons les sommes de G‡l
$$S_\alpha( \M):=\sum_{m,n\in\M} {(m,n)^\alpha\over [m,n]^\alpha}\qquad (\alpha>0).\eqdef{defSa}$$
Ces expressions faisant appara"tre des plus grands communs diviseurs ont ŽtŽ introduites par  Erd\H os.  G‡l \citer{G49} a rŽsolu la conjecture d'Erd\H os correspondante dans le cas  $\alpha=1$.
De nombreux articles rŽcents 
\citer{ABS15}, \citer{BHS16}, \citer{BS15}, \citer{BS17},
\citer{LR17} 
concernent le comportement asymptotique de la quantitŽ
$$\Gamma_\alpha(N):=\sup_{|\M|= N}{S_\alpha( \M)\over |\M|},\note{En fait, dans les travaux citŽs, le supremum est pris sur $|\M|\leqslant N$. La dŽfinition adoptŽe ici permet un gain de prŽcision.}\eqdef{defGa}$$
elles-mme liŽe aux majorations de certains polyn™mes de Dirichlet et ˆ celles de maximums localisŽs de la fonction zta de Riemann sur la droite verticale d'abscisse $\alpha$.  \par 
Soit $\N_1$ l'ensemble des entiers sans facteur carrŽ. Nous notons $\Gamma^*_\alpha(N)$ la quantitŽ analogue ˆ $\Gamma_\alpha(N)$ obtenue en imposant $\M\subset \N_1$. \par 
Nous avons essentiellement restreint la prŽsente Žtude au cas $\alpha=\dm$, qui se rŽvle tre l'un des plus intŽressants --- voir le survol \citer{SS16}.  Nous posons en consŽquence $S(\M)=S_{1/2}(\M)$, $\Gamma (N)=\Gamma_{1/2}(N)$ et~$\Gamma^* (N)=\Gamma^*_{1/2}(N)$.
\par 
Soit $$\L(x):=\exp\Bigg\{\sqrt{{\log x \log_3x\over \log_2 x}}\Bigg\}\qquad (x>16),$$
o, ici et dans la suite, nous notons $\log_k$ la $k$-ime itŽrŽe de la fonction logarithme.
Dans~\citer{BS15}, Bondarenko et Seip ont Žtabli l'existence  d'une constante $A>0$ telle que l'on ait
$$\Gamma(N)\leqslant \L(N)^A\qquad (N>N_0(A)).\eqdef{majGamma}$$
Il est indiquŽ dans le mme travail que la valeur $A=7$ est admissible. Il y est Žgalement Žtabli que  \eqref{majGamma}  n'est pas valable pour~$A<1$ (\citer{BS17}, thŽorme 2).

\medskip
Nous nous proposons ici de prŽciser ces rŽsultats.
\Propt{th}{Lorsque $N$ tend vers l'infini, nous avons
$$\Gamma (N)= \L(N)^{2\sqrt{2}+o(1)} .\eqdef{minGamma}$$ }

 \rem Notre approche fournit Žgalement l'estimation
$$\Gamma^* (N)= \L(N)^{2 +o(1)}.\eqdef{maj^*Gamma}$$

 Notons
$ \no{\bfc}_2$  la norme  quadratique  d'une suite complexe $\bfc:=\{c_n\}_{n=1}^{\infty}$. Posons
$$Q(\M):=\sup_{\di{\bfc\in \CC^N}{\no{\bfc}_2=1}}
\abs{ \sum_{m,n\in \M}c_m\ov{c_{n}}\sqrt {(m,n) \over [m,n]}}.\eqdef{defQM}$$
 En vertu du thŽorme 5  de~\citer{ABS15}, nous avons
$$
\Gamma(N)\leqslant \sup_{|\M|=N}Q(\M)\leqslant (\e^2+1)(\log N+2)\Gamma^+(N),\eqdef{encth5}
$$o $\Gamma^+ (N):=\max_{n\leqslant N} \Gamma(n).$
Alors que l'inŽgalitŽ de gauche est triviale, celle de droite est remarquable : elle exprime que  certains des  vecteurs propres de la forme quadratique associŽe ˆ une somme de G‡l maximale sont proches  de la droite vectorielle engendrŽe par le vecteur $(1,\ldots,1)$. \par 
Il est ˆ noter que, quitte ˆ remplacer $\Gamma(N)$ par $\Gamma^*(N)$, l'encadrement \eqref{encth5} persiste lorsque le supremum est restreint aux ensembles $\M\subset \N_1$.
\par  
L'estimation suivante dŽcoule immŽdiatement de \eqref{encth5}.
\Propt{cor1}{Lorsque $N$ tend vers l'infini, nous avons
$$\sup_{|\M|=N}Q(\M)=\L(N)^{2\sqrt{2}+o(1)} . $$ }

\par \medskip
 La structure multiplicative des sommes de G‡l est plus simple lorsque l'ensemble $\M$ est choisi comme
 l'ensemble $\T_D $ de tous les diviseurs d'un entier $D$. Ainsi que l'on peut s'y attendre, les analogues des maximums prŽcŽdents subissent une rŽduction significative.  Distinguons deux cas selon que $D$ est ou non sans facteur carrŽ: 
$$\Gamma'(N):=\sup_{\tau(D)\leqslant N}{S(\T_D)\over N},\qquad \Gamma''(N):=\sup_{\tau(D)\leqslant N}{\mu(D)^2S(\T_D)\over N}\cdot$$
Soit $B$ la constante dŽfinie par $$B:=4\sqrt{\sum_{k\geqslant 1} {1\over k^2(1+k)^2\log (1+1/k) } }\cdot\eqdef{defA}$$
\Propt{maj}{Lorsque $N$ tend vers l'infini, nous avons
$$\leqalignno{\log \Gamma'(N)&=    B	\sqrt{\log N\over \log_2N}\Big\{1+O\Big({\log_3N\over \log_2N}\Big)\Big\},&\eqdef{G'} \cr
\log \Gamma''(N)&= 
{2\over \sqrt{\log 2}}\sqrt{\log N\over \log_2N}\Big\{1+O\Big({\log_3N\over \log_2N}\Big)\Big\} .&\eqdef{G""}\cr}$$  }

\medskip

\rems (i) On voit que $B>2/\sqrt{\log 2}$ en considŽrant le premier terme de la sŽrie de \eqref{defA}. On a en fait
$B\approx2,78422$ alors que $2/\sqrt{\log 2}\approx2,40224$\par \smallskip
(ii) Une minoration du type  $\log \Gamma''(N)\geqslant c\sqrt{(\log N)/\log_2N}$ est Žtablie dans~\citer{ABS15}.
\bigskip\goodbreak
\paradeux{Maximums localisŽs de la fonction zta}
 Dans le prolongement des travaux  de Montgomery \citer{M77} et Balasubramanian--Ramachandra \citer{BR77}, Bondarenko et Seip ont rŽcemment obtenu --- voir notamment  \citer{BS17}, \citer{BS17b} ---  des bornes infŽrieures pour la quantitŽ 
$$Z_\beta(T):=\max_{T^\beta\leqslant \tau\leqslant T} \big| \zeta(\dm+i\tau)\big|\qquad (0\leqslant \beta<1,\,T\geqslant 1)\eqdef{defZT}$$
 en minorant  la norme de la forme quadratique associŽe ˆ la sous-somme $\gotS(\M)$.          
Leur approche repose sur la mŽthode de rŽsonance, dŽveloppŽe  indŽpendam\-ment,  dans ce cadre, par Soundararajan \citer{S08} et Hilberdink \citer{Hi09}. Elle est adaptable ˆ l'Žtude des maximums d'autres  sŽries de Dirichlet --- voir la publication rŽcente \citer{BS17b}.
 
En adaptant la mŽthode de \citer{BS17}, nous dŽduisons de nos estimations le rŽsultat suivant.
\Propt{thzeta}{
Soit $\beta\in [0,1[$ et $c$ une constante vŽrifiant $0<c< \sqrt{2(1-\beta)} .$ Lorsque $T$ est suffisamment grand, nous avons
$$Z_\beta(T)\geqslant \L(T)^{c}.$$}

\rems (i) Dans \citer{BS17}, il est montrŽ que $c< \sqrt{\min\{ \dm, 1-\beta\}}$ est une condition admissible. Dans \citer{BS17c}, ce rŽsultat est Žtendu ˆ  $c< \sqrt{  1-\beta }$ ce qui reprŽsente  une amŽlioration d'un facteur $\sqrt{2(1-\beta)}$  lorsque $\beta\leqslant \dm$. Le \ref{thzeta} amŽliore donc les exposants obtenus par Bondarenko et Seip d'un nouveau facteur $\sqrt{2}$ pour tout $\beta\in [0,1[$. \par 
(ii) Notre approche diffre sensiblement de celle de \citer{BS17c}: ˆ l'instar des mŽthodes dŽveloppŽes dans \citer{Hi09} et \citer{Ai16}, nous relions le carrŽ du maximum ˆ la somme de G‡l complte $S(\M)$.
L'introduction du carrŽ explique que le gain obtenu n'est que la moitiŽ de celui qui a ŽtŽ obtenu pour la somme de G‡l soit $\dm 2\sqrt{2}=\sqrt{2}.$  Voir Žgalement la remarque ˆ la fin du \S\thinspace\ref{prgvZ}. 
\medskip
\medskip
\paradeux{Valeurs maximales de $|L(\dm,\chi)|$}\smallskip
ConsidŽrons les fonctions $L$ associŽes aux caractres de Dirichlet $\chi$, soit  $$L(s,\chi):=\sum_{n\geqslant 1} {\chi (n)\over n^s}\qquad (\sigma:=\re (s)>1).$$
Notons
$$X_q^+:=\{ \chi (\mod q)   :  \chi(-1)=1\},\qquad X_q^-:=\{ \chi (\mod q):  \chi(-1)=-1\}.$$
Le \ref{th} peut tre combinŽ ˆ la mŽthode de rŽsonance pour obtenir une minoration des quantitŽs
$$  L_q^+:=\max_{\di{\chi \in X_q^+}{\chi\neq \chi_0}}\vbs{L(\dm,\chi)}\qquad (q\geqslant 3),$$
o $\chi_0$ dŽsigne le caractre principal modulo $q$.
Nous nous limitons ici au cas o $q$ est un \np. La restriction aux caractres pairs provient ici de la forme particulire de la relation d'orthogonalitŽ des caractres de paritŽ donnŽe (\cf~formule \eqref{orthprim} {\it infra}) et de la nature des coefficients obtenus lors du dŽveloppement des sommes pondŽrŽes de valeurs $|L(\dm,\chi)|^2$: \cf~\eqref{V2+}~{\it infra}.
\Propt{thL}{Lorsque le \np\ $q$ tend vers l'infini, nous avons
$$L_q^+\geqslant \L(q)^{1+o(1)}.$$}
\par 
 Des estimations de mme nature ont ŽtŽ abondamment considŽrŽe dans la littŽrature. Citons-en quelques unes. 
\par 
Dans \citer{S08}, Soundararajan a Žtabli que
$$\max_{x<|d|\leqslant 2x}\log L(\dm,\chi_d)\geqslant (1+o(1))\sqrt{\log x\over 5\log_2x}\qquad (x\to\infty),$$
o $d$ parcourt la suite des discriminants fondamentaux et $\chi_d$ dŽsigne le caractre rŽel associŽ ˆ $d$ (nous avons rectifiŽ une faute de frappe apparaissant dans l'ŽnoncŽ correspondant de \citer{S08}).
\par \goodbreak
Dans \citer{H16}, Hough montre que, pour tout \np\ $q$  suffisamment grand, $$\max_{\di{\chi\mod q}{\chi\neq\chi_0}}\log \vbs{L(\dm,\chi)}\geqslant {1\over 4} \sqrt{\log q\over 2\log_2q} $$ avec la prŽcision supplŽmentaire que l'argument de $L(\dm,\chi)$ peut tre choisi proche de tout $\theta\in [0,2\pi]$ fixŽ ˆ l'avance. On se reportera ˆ \citer{H16} pour un ŽnoncŽ prŽcis.
\par \goodbreak
 Des minorations de $|L(\sigma,\chi)|$ ont ŽtŽ Žtablies rŽcemment dans \citer{AKMP18}  pour $\sigma\in ]\dm,1 [$. Dans ce travail, il est notamment prouvŽ que 
$$\max_{\di{\chi\mod q}{\chi\neq\chi_0}}\log \vbs{L(\sigma,\chi)}\gg_\sigma  (\log q)^{1-\sigma}(\log_2q)^\sigma\qquad (q\geqslant 3,\,\dm<\sigma< 1). $$ Cette minoration avait prŽcŽdemment ŽtŽ Žtablie  dans \citer{L11} sous l'hypothse  de Riemann gŽnŽralisŽe. 
\bigskip

\paradeux{Grandes valeurs de sommes de caractres}
Posons, pour tout caractre de Dirichlet $\chi$,
$$S(x,\chi):=\sum_{n\leqslant x} \chi(n).$$
Dans \citer{H13}, amŽliorant des estimations de \citer{GS01}, Hough Žtablit une minoration de la quantitŽ 
$$\Delta(x,q):=\max_{\di{\chi\neq \chi_0}{\chi\mod q}}\vbs{S(x,\chi)},$$
lorsque $q$ est un nombre premier. Notant $x=q^\vartheta$, cette estimation est valide dans le domaine $$4\sqrt{\log_2q\over \log q}\log_3q\leqslant \vartheta\leqslant 1-4\sqrt{\log_2q\over \log q}\log_3q.\eqdef{domH}$$
\par 
Nous pouvons dŽduire rapidement du \ref{th} le rŽsultat suivant, valide sans restriction sur le module $q$. Ici et dans la suite, nous notons $\omega(q)$ le nombre des facteurs premiers distincts d'un entier $q$.
\Propt{hough}{Soit $\varepsilon>0$. Sous la condition $\e^{(\log q)^{1/2+\varepsilon}}\leqslant x\leqslant q/\e^{(1+\varepsilon)\omega(q)}$, nous avons
$$\Delta(x,q) \gg  \sqrt{x}\L(q/x )^{ \sqrt{2}+o(1)}\qquad (q\to\infty).\eqdef{H+}
$$}

\rem On a en toute circonstance $q/\e^{(1+\varepsilon)\omega(q)}\geqslant q^{1-\{1+ \varepsilon+o(1)\}/\log_2q}$ lorsque $q\to\infty$. 
\medskip
Dans son domaine de validitŽ, cette minoration amŽliore celle de Hough (\citer{H13}, thŽorme 3.1) d'un facteur  $\sqrt{2\log_3(q/x)}$ dans l'exposant. Cependant, lorsque $q$ est premier, le domaine \eqref{domH}, plus grand que celui de \eqref{H+},  englobe celui du changement de phase observŽ autour de  $$\log x=\sqrt{\log q\log_2q}.$$  Nous renvoyons le lecteur ˆ \citer{H13} pour plus de dŽtails.
Lorsque $\vartheta>1-C/\log_3q $, la minoration donnŽe dans \citer{H13} est meilleure que celle de \eqref{H+}. 
Dans \citer{H11}, Hough considre Žgalement le cas de modules $q$ composŽs, sous la condition $\vartheta\leqslant 1-\varepsilon$.

\medskip

\paraun{Preuve du \ref{th}: minoration}

\paradeux{RŽduction du problme} \drefun{reduc} 

Le rŽsultat suivant permet de rŽduire la minoration contenue dans \eqref{majGamma} et \eqref{maj^*Gamma} aux cas de  supremums pour $|\M|\leqslant N$.  
\Propl{prelim}{
Soient $N\in\N$ et $\M$ un ensemble d'entiers de cardinal $\leqslant N$. Alors, il existe  un ensemble d'entiers $\M'$ de cardinal $N $ tel que
$${S(\M')\over |\M'|}\geqslant {S(\M)\over 2|\M |}\cdot$$
De plus, si $\M\subset \N_1$, on peut choisir $\M'\subset\N_1$. }
\medskip
\nid
La premire Žtape consiste ˆ montrer qu'il   existe un ensemble $\M''$ de cardinal $N'' \in ]\dm N,N]$   tel que
$${S(\M'')\over |\M''|}\geqslant {S(\M)\over |\M|}\cdot$$ 
Si $N/2<|\M|\leqslant N$, il n'y a rien ˆ dŽmontrer. Si $|\M|\leqslant N/2$, il existe un entier $k\geqslant 1$ tel que $N/2\leqslant 2^k|\M|\leqslant N$. ConsidŽrons  alors un ensemble $\P$ de $k$ nombres premiers distincts ne divisant aucun ŽlŽment de $\M$. Posant $D:=\prod_{p\in\P}p$, nous dŽfinissons $$\M'':=\big\{ dm : d\mid D,\, m \in \M \big\}.$$
Nous avons d'une part $|\M''|= 2^k |\M|\in  [N/2,N]$ et, d'autre part,
$$S(\M'' )=\sum_{d,d'\mid D}{(d,d')\over \sqrt{dd'}}S(\M' )\geqslant  2^kS(\M)
= |\M''|{S(\M)\over |\M|},$$ 
 o l'inŽgalitŽ est obtenue en restreignant la sommation ˆ $d=d'$.
\par 
Pour achever la preuve, il suffit de complŽter  l'ensemble $\M''$ en un ensemble $\M'$ de cardinal $N $ et d'observer que  
$ S(\M')\geqslant S(\M'')$.
\qed
\medskip
\paradeux{Construction d'un ensemble $\M$}
 Ce paragraphe est consacrŽ ˆ la construction d'ensembles $\M$ pour lesquels  la somme  $S(\M)$  atteint de grandes valeurs. Dans \citer{BS17} les ensembles construits Žtaient constituŽs d'entiers sans facteur carrŽ.  Notre amŽlioration est  notamment due  au fait que nous nous affranchissons de cette restriction. Par ailleurs,  nous minorons directement $S(\M)$ alors que la mŽthode employŽe dans  \citer{BS17} nŽcessite de considŽrer la forme quadratique associŽe.
\par 
Dans toute la suite, nous notons $\PP$ l'ensemble de tous les \nps\ et nous considŽrons un entier $N$ arbitrairement grand. Soient $u\in ]1,\e]$, $a\in ]1,+\infty[$ et $\gamma\in ]0,1[$ trois paramtres  bornŽs  dont les valeurs seront choisies ultŽrieurement.
Introduisons les intervalles
$$I_k:=\Big]u^k\log N\log_2N, u^{k+1}\log N\log_2N\Big]\qquad \big(1\leqslant k\leqslant (\log_2N)^\gamma\big),$$  
de sorte que
$$\eqalign{P_k&:=|I_k\cap\PP|=\pi(u^{k+1}\log N\log_2N)-\pi(u^{k}\log N\log_2N)
\leqslant  u^{k +1} \log N  .\cr}\eqdef{estPk}$$  
De plus
$$ P_k = u^{k }(u-1)\log N \Big\{1+O\Big({k+\log_3N\over \log_2N}\Big)\Big\}.$$
Posons encore
 $$J_k:=2\fl{a\log N\over 2k^2\log_3N}\qquad \big(1\leqslant k\leqslant (\log_2N)^\gamma\big)\eqdef{defJk}$$ et choisissons $a>1$  tel que  $a\gamma<1/\log u$.
\par 

Notant $N_k:=\prod_{p\in I_k}p$, nous considŽrons les ensembles
$$\eqalign{ 
\M_k  &:=\Big\{ m  \,:\, m={\ell\over q}N_k,\,\omega(\ell)\leqslant \dm J_k,\,\omega(q)\leqslant \dm J_k, \,\ell q\mid N_k\Big\} \cr}$$   
et nous choisissons
$$\M :=\Big\{ m=\prod_{1\leqslant  k\leqslant (\log_2N)^\gamma} m_k\,:\quad m_k\in \M_k \quad  \big(1\leqslant k\leqslant (\log_2N)^\gamma\big)\Big\}.\eqdef{defM}$$

\par 
\Propl{lemM}{Sous la condition $a<1/(\gamma \log u)$, nous avons  $ |\M |\leqslant N.$
}\medskip
\nid Nous avons  
$$|\M|=\prod_{1\leqslant k\leqslant (\log_2N)^\gamma} |\M_k|.\eqdef{cardM=prod}$$
avec
$\dsp|\M_k|=\sum_{\di{0\leqslant j\leqslant J_k/2}{0\leqslant h\leqslant J_k/2}}
\Big({P_k\atop j}\Big)\Big({P_k-j\atop h}\Big). 
$ 

Les relations
$${P_k\choose j-1}\leqslant \dm {P_k\choose j},\quad 
{P_k-j\choose h-1}\leqslant \dm {P_k-j\choose h}\qquad \big( 1\leqslant j, h\leqslant \dm J_k\leqslant \ft 14 P_k\big)$$
et $$
{P_k\choose j}{P_k-j\choose J_k/2}={P_k \choose J_k/2}{P_k-J_k/2\choose j}\qquad \big(0\leqslant j\leqslant J_k/2\big)$$
fournissent 
$$ \eqalign{{P_k\choose J_k/2}{P_k-J_k\choose J_k/2}\leqslant |\M_k|
&\leqslant  2\sum_{{0\leqslant j\leqslant J_k/2} }
{P_k\choose j}{P_k-j\choose J_k/2}\cr
&\leqslant  4{P_k\choose J_k/2}{P_k-J_k/2\choose J_k/2}.\cr}\eqdef{encadMk} $$

Gr‰ce ˆ l'inŽgalitŽ  $${m\choose n}\leqslant  \e^n\Big({m\over n}\Big)^n  \qquad (1\leqslant n\leqslant m) $$ 
qui dŽcoule par exemple de la formule d'Euler-Maclaurin, et ˆ \eqref{estPk}, nous obtenons, pour $N$ assez grand,
$$|\M_k|\leqslant 4\Big({2\e P_k\over J_k}\Big)^{J_k} \leqslant  \Big(\e k^2u^{k+1}\log _3N\Big)^{J_k} .$$
En reportant dans \eqref{cardM=prod}, cela implique
$$|\M|\leqslant N^{a\gamma \log  u+o(1)}.$$   
\vskip-5mm\qed
\medskip
\paradeux{ComplŽtion de l'argument}
Il reste ˆ minorer le rapport $ S(\M)/|\M|$ lorsque $\M$ est dŽfini par \eqref{defM} et les paramtres $a,$ $u$, $\gamma$ sont choisis comme indiquŽ plus haut.
Nous avons $$S (\M)=\prod_{1\leqslant k\leqslant (\log_2N)^\gamma} S(\M_k).\eqdef{SM=prodSMk}$$
\par 
Lorsque les entiers
$$m:= {\ell\over  q}N_k, \qquad m':= {\ell'\over q'}N_k$$
sont dans $\M_k$, nous avons
$$(m,m')={N_k\over [q,q']}\Big({\ell q'\over (q,q')},{\ell'q\over (q,q')}\Big)={N_k(\ell,\ell')\over [q,q']},$$
et donc
$${m\over (m,m')}={\ell\over (\ell,\ell')}{q'\over (q,q')},\qquad {m'\over (m,m')}={\ell'\over (\ell,\ell')}{q\over (q,q')}\cdot$$
Ainsi
$$\eqalign{
S (\M_k )&=\sum_{\di{\ell,\ell'\mid N_k}{\omega(\ell ),\omega(\ell ')\leqslant J_k/2}} {(\ell,\ell') \over \sqrt{\ell \ell'}}
\sum_{\tri{q,q'\mid N_k}{(q,\ell)=(q',\ell')=1}{\omega(q ),\omega(q')\leqslant  J_k/2}}  {(q,q') \over \sqrt{qq'}} 
\cdot\cr}$$
DŽsignons par
$S_k=S_k(\ell,\ell')$ la somme intŽrieure.
L'identitŽ
$n=\sum_{d\mid n}\phi(d)$ permet d'Žcrire
$$\eqalign{S_k&=\sum_{\tri{d_1\mid N_k}{\omega(d_1)\leqslant J_k/2}{(d_1,\ell\ell')=1}}  {\phi(d_1)\over d_1} 
\sum_{\tri{n_1,n'_1\mid N_k }{(n_1,\ell d_1)=(n'_1,\ell'd_1)=1}{\omega(n_1 ),\omega(n'_1)\leqslant  J_k/2-\omega(d_1)}}  {1\over \sqrt{n_1n'_1}} \cr
&=\sum_{\tri{d_1\mid N_k}{\omega(d_1)\leqslant J_k/2}{(d_1,\ell\ell')=1}}  {\phi(d_1)\over d_1} \sigma(\dm J_k-\omega(d_1),\ell d_1)\sigma(\dm J_k-\omega(d_1),\ell' d_1) 
\cr}$$
o l'on a posŽ
$$\sigma(R,r ):=
\sum_{\tri{ n\mid N_k }{(n,r)=  1}{\omega(n ) \leqslant  R}}  {1\over \sqrt{n }}\qquad (r\geqslant 1,\,R\geqslant 1).$$
\par 
Comme
$${\phi(d)\over d}\gg  {\exp\Big\{-\sum_{p\in I_k}{1\over p}\Big\}}\qquad (d\mid N_k),$$  
il suit
$$S_k\gg\sum_{\tri{d_1\mid N_k}{\omega(d_1)\leqslant J_k/2}{(d_1,\ell\ell')=1}} \sigma(\dm J_k -\omega(d_1),\ell d_1)\sigma(\dm J_k-\omega(d_1),\ell' d_1) . \eqdef{est}$$
\par 
En utilisant la mme identitŽ pour exprimer  $(\ell,\ell')$, nous obtenons
$$\eqalign{
S (\M_k )\gg\hskip-3mm
\sum_{\tri{d_1,d_2\mid N_k}{\omega(d_j)\leqslant J_k/2}{(d_1,d_2)=1} }   \!\!
\sum_{\tri{n_2,n'_2\mid N_k }{(n_2 n'_2, d_1d_2)=1}{\omega(n_2 ),\omega(n'_2)\leqslant J_k/2-\omega(d_2)}} \hskip-5mm {\sigma(\dm J_k -\omega(d_1),n_2d_2 d_1)\sigma(\dm J_k-\omega(d_1),n_2'd_2 d_1) \over \sqrt{n_2n'_2}}\cdot\cr}$$
Soit
 $$j_k:=\fl{{\alpha \over k}\sqrt{{\log N\over \log_2N\log_3N}}}$$ o $\alpha $  est un paramtre bornŽ  ˆ optimiser.   Restreignons la somme extŽrieure aux couples $(d_1,d_2)$ tels que $\omega(d_j)\leqslant \dm J_k-j_k$  et la somme intŽrieure aux  entiers $n_2$, $n_2'$ tels que $\omega(n_2)=\omega(n_2')=j_k$. \par 
Pour $n=n_2$ ou $n=n'_2$, nous avons
$$ \sigma(\dm J_k -\omega(d_1),nd_2 d_1)\geqslant \sigma(j_k,nd_2 d_1)\gg {1\over j_k !}\Bigg(\sum_{\di{p\in I_k}{p\nmid  nd_2 d_1}}{1\over \sqrt{p}}\Bigg)^{j_k}.\eqdef{minsigma}$$
Or
$$\sum_{p\in I_k}{1\over \sqrt{p}}=2u^{k/2}(\sqrt{u}-1)\sqrt{{\log N\over \log_2N}}\Big\{1+O\Big({k+\log_3N\over \log_2N}\Big)\Big\}$$
tandis que, lorsque  $d_1d_2n\mid N_k$  et $\omega(d_1)+\omega(d_2n)\leqslant J_k$,
$$ \sum_{p\mid d_1d_2n}{1\over \sqrt{p}}\ll {J_k\over u^{k/2}\sqrt{\log N\log_2N}}\ll {1 \over k^2u^{k/2} \log_3N}\sqrt{{\log N\over \log_2N}}\cdot$$
Nous pouvons donc Žcrire 
$$\sum_{\di{p\in I_k}{p\,\nmid\, d_1d_2n}}{1\over \sqrt{p}}=2u^{k/2}(\sqrt{u}-1)\sqrt{{\log N\over \log_2N}}\Big\{1+O\Big({1\over \log_3N}\Big)\Big\}.\eqdef{estsumpnmidd}$$
En utilisant la formule de Stirling sous la forme 
$ j_k!= \exp\big\{ j_k\log j_k-j_k+ O(\log j_k)\big\}  $ et en reportant dans  \eqref{minsigma},
nous obtenons donc
 $$\sigma(j_k,nd_2 d_1)\geqslant T^{j_k}\e^{o(j_k)}\eqdef{minsig}$$
avec
$$\eqalign{T&:={2\e(\sqrt{u}-1)u^{k/2}\over j_k}\sqrt{\log N\over \log_2N}\cr&={2k\e(\sqrt{u}-1)u^{k/2}\over \alpha}\sqrt{\log_3N}\bigg\{1+O\bigg({\sqrt{\log_2N\log_3N}\over \sqrt{\log N}}\bigg)\bigg\},\cr}\eqdef{estT}$$   
puis
$$\eqalign{{
S (\M_k ) }\gg  T^{{2j_k}}\e^{o(j_k)}
\sum_{\tri{d_1,d_2\mid N_k}{\omega(d_j)\leqslant J_k/2-j_k}{(d_1,d_2)=1} }   \sum_{\tri{n_2,n'_2\mid N_k }{(n_2 n'_2, d_1d_2)=1}{\omega(n_2 )=\omega(n'_2)=j_k}} {1\over \sqrt{n_2n'_2}}\cdot\cr}$$
Par un calcul parallle ˆ celui qui fournit \eqref{minsig}, l'estimation  \eqref{estsumpnmidd} permet de montrer que la somme intŽrieure vaut
$$\Bigg(\sum_{\tri{n \mid N_k }{(n  , d_1d_2)=1}{\omega(n ) =j_k}} {1\over \sqrt{n }}\Bigg)^2\gg  T^{2j_k}\e^{o(j_k)}.$$
 \par 
 Nous avons ainsi Žtabli que 
 $$S(\M_k)\gg T^{4j_k}\e^{o(j_k)}V_k\eqdef{est1SMk}$$ 
avec
 $$\eqalign{V_k:=&\sum_{\tri{d_1,d_2\mid N_k}{\omega(d_j)\leqslant J_k/2-j_k}{(d_1,d_2)=1} }1
=
\sum_{\di{0\leqslant j\leqslant J_k/2-j_k}{0\leqslant h\leqslant J_k/2-j_k}}
{P_k\choose j}{P_k-j\choose h}\geqslant 
{P_k\choose J_k/2-j_k}{P_k-J_k/2+j_k\choose J_k/2-j_k}\cr
&={P_k\choose J_k/2 }{P_k-J_k/2 \choose J_k/2 }{(J_k/2)!^2(P_k-J_k)!\over (J_k/2-j_k)!^2(P_k-J_k+2j_k)!}\cr
&\geqslant \ft14|\M_k| \Big({J_k\over 2P_k}\Big)^{2j_k}\geqslant |\M_k| \bigg({a+o(1)\over 2k^2u^k(u-1)\log_3N}\bigg)^{2j_k} ,\cr}\eqdef{estVk}$$
en vertu de \eqref{encadMk}. 
 \par 
 Gr‰ce ˆ \eqref{estT}, \eqref{est1SMk} et \eqref{estVk}, nous pouvons Žcrire
$${
S (\M_k )\over |\M_k |}\gg \bigg({h\over \alpha^2}\bigg)^{2j_k} \e^{o(j_k)} $$
o l'on  a posŽ $h:={2\e^2a(\sqrt{u}-1)/(\sqrt{u}+1}).$
En reportant \eqref{SM=prodSMk}, nous obtenons
$${S(\M)\over |\M|}\geqslant \L(N)^{\beta+o(1)}\eqdef{minSMb}$$
avec $\beta:=2\gamma \alpha\log (h/\alpha^2).$ Le choix optimal $\alpha:=\sqrt{h}/\e$ fournit
$$\beta ={4\gamma\sqrt{h}\over \e}={4\gamma}\sqrt{{2a(\sqrt{u}-1)}\over \sqrt{u}+1}\cdot$$
En choisissant  $a\gamma\log u$ proche de 1 puis $u$ proche de 1, nous constatons que $\beta$ peut tre pris arbitrairement proche de $2\sqrt{2}$. 
\par 
Cela achve la dŽmonstration de la minoration du \ref{th}.

\paraun{Preuve  du \ref{th}:  majoration}
 
\paradeux{Majoration de sommes pondŽrŽes}\drefdeux{majsp}

Nous adaptons ici la mŽthode dŽveloppŽe par Bondarenko et Seip dans \citer{BS15}, qui repose sur les rŽsultats de \citer{ABS15}.
ConsidŽrons les sommes de G‡l pondŽrŽes
$$S(\M;g):=\sum_{m,n\in\M}  g\bigg({[m,n]\over (m,n)}\bigg)\eqdef{galpond}$$
o $g$ est une fonction sous-multiplicative,\note{Autrement dit telle que $g(mn)\leqslant g(m)g(n)$ pour tous entiers $m,n$ premiers entre eux.} de sorte que
$$\eqalign{S(\M;g)\leqslant S^+(\M;g)&:= \sum_{m,n\in\M}g\Big({m\over (m,n)}\Big)g\Big({n\over (m,n)}\Big)=\sum_{\ell}\sum_{\di{m,n\in\M}{[m,n]=\ell}}g\Big({\ell\over m}\Big) g\Big({\ell\over n}\Big) \cdot\cr}
\eqdef{majgenSf}$$
 D'aprs \citer{BS15}, nous avons
$$\Gamma(N)\leqslant 4\sup_{\di{ |\M|\leqslant N}{\M\subset \N_1}} {S\big(\M;2^\omega g_0\big)\over|\M| } ,\eqdef{majC}$$
o, ici et dans la suite, $\omega$ dŽsigne la fonction nombre de facteurs premiers, comptŽs sans multiplicitŽ, et  $g_0(n):=1/\sqrt{n}$ $(n\geqslant 1$).
\par 

Comme dans \citer{BS15},  nous convenons qu'un ensemble d'entiers $\M$ est dit {\it divisoriellement clos} s'il contient tous les diviseurs de chacun de ses ŽlŽments, et qu'il est dit {\it complet} si,  pour chaque ŽlŽment $m=\prod_{j=1}^r p_j^{\nu_j}$ de $\M$ o les $p_j$ sont distincts l'entier $n:=\prod_{j=1}^r {q}_j^{\nu_j}$ est Žgalement dans $\M$  ds que les ${q}_j $ sont distincts et satisfont $q_j\leqslant p_j$. 
\par 
Nous dirons qu'un ensemble fini d'entiers est  {\it strict} s'il est divisoriellement clos et complet.
D'aprs les lemmes~1 et 2 de \citer{BS15}, pour chaque $\alpha$ de $]0,1[$, le supremum  $\Gamma^{*}_\alpha(N)$  est atteint pour un ensemble strict $\M$. 
\par 
Le lemme-clef de notre dŽmonstration dŽcoule d'un argument de complŽtude, essentiellement identique au lemme 3 de \citer{BS15}.  
\par 
Dans toute la suite, nous notons $\{p_j\}_{j\geqslant 1}$ la suite croissante des nombres premiers et posons $y:=p_{\fl{\log N/\log 2}}$.

\Propl{lem3}{Soient  $\M$ un ensemble strict de $N$ entiers distincts et $n\in\M$. Si $$y<p_{j_1} <p_{j_2} <\ldots<p_{j_\nu}$$ sont les facteurs premiers de $n$ excŽdant $y$, alors
$$\sum_{1\leqslant h\leqslant \nu} \log (j_h/2 h)  \leqslant  \log N.\eqdef{bs3}$$ 
De plus, $\nu=\nu_n\leqslant (\log N)/\log 2$.
}
\medskip
 \nid Comme indiquŽ plus haut, il s'agit d'une variante du lemme 3  de \citer{BS15}. Nous donnons les dŽtails pour la commoditŽ du lecteur. \par 
Il existe au moins $R:=\prod_{1\leqslant h\leqslant \nu}(j_h-h+1)/\nu!$ entiers de la forme $\prod_{1\leqslant h\leqslant \nu}p_{s_h}$ avec $\nu<s_h\leqslant j_h$ et tels que la suite $\{s_h\}_{h=1}^{\nu}$ soit strictement croissante. Comme $\M$ est complet, nous avons $R\leqslant N$. Comme de plus $\M$ est divisoriellement clos, nous avons Žgalement $2^\nu\leqslant N$. En observant que $(h-1)/j_h\leqslant (h-1)/(\nu+h)\leqslant \dm$ pour $1\leqslant h\leqslant \nu$, il suit
$${1\over 2^\nu\nu!}\prod_{1\leqslant h\leqslant \nu}j_h\leqslant R\leqslant N.$$
On obtient bien \eqref{bs3} en prenant les logarithmes.\qed 
\par\goodbreak 
\medskip
D\'{e}signons par $g_1$ la fonction arithmŽtique multiplicative dŽfinie par $$g_1(n):={\mu(n)^2\over \prod_{p\mid n}(\sqrt{p}-1)}\cdot $$ 
 \Propl{lemmessfc}{Soient $\M\subset \N_1$ un ensemble de $N$ entiers  sans facteur carrŽ et $C\geqslant 1$. Alors
 $$S (\M;C^\omega g_1)\leqslant N \L(N)^{2C+o(1)}.\eqdef{ssfc}$$}
\par 
\rems (i)
 Dans \citer{BS15}  --- \cf\ Žquation (12) ---, cette majoration est ŽnoncŽe avec un exposant $2\sqrt{6C}$  mais il semble que la dŽmonstration contienne une erreur et ne fournisse que $2\sqrt{6}C$.  Par consŽquent, il faut modifier la majoration finale $\Gamma(N)\leqslant \L(N)^{A+o(1)}$ obtenue dans \citer{BS15}  en remplaant $A=2\sqrt{12}<7$ par $A=4\sqrt{6}<10$.
\par \smallskip
(ii)
Notre approche permet d'employer la majoration \eqref{ssfc} avec $C=\sqrt{2}$ au lieu de $C=2$ comme dans \citer{BS17}, ce qui fournit un gain de prŽcision. \smallskip
\par  (iii) La mŽthode du \S\ref{prgreduc} {\it infra} fournit, pour chaque $N\geqslant 1$, l'existence d'un ensemble  $\M_N\subset \N_1$ de cardinal $N$ et tel que 
 $$S \big(\M_N;C^\omega g_1\big)\geqslant N \L(N)^{2C+o(1)}. $$ Nous omettons les dŽtails.
 
\medskip

\noi \it DŽmonstration du \ref{lemmessfc}. \rm
En vertu des rŽsultats de \citer{ABS15} et \citer{BS15} rappelŽs plus haut, nous pouvons supposer que $\M$ est strict. 
Nous dŽduisons deux propriŽtŽs de cette hypothse.\par  
Lorsque les $p_{j_h}$ sont les $\nu$ facteurs premiers $>y$ d'un entier $m\in \M$, nous avons donc,  d'aprs le \ref{lem3},
$$\sum_{1\leqslant h\leqslant \nu} \big\{\log (j_h)- \log (2\log N/\log 2 ) \big\} \leqslant  \log N,$$
et $\nu\leqslant (\log N)/\log 2$. Comme  $\log j=\log (p_j/\log p_j)+O(1),$ 
il existe une constante $B>0$ telle que
$$\sum_{\di{p\mid m}{p>y}}  \log \Big({p\over B\log N\log p}\Big) \leqslant  \log N\qquad (m\in\M).\eqdef{factmy}$$
\par 
Par ailleurs, si $n\in \M$ possde un facteur premier plus grand que $p_{N}$, ce nombre premier est dans~$\M$, car $\M$ est divisoriellement clos. Cela implique que tous les nombres premiers $p_j$ avec $j\leqslant N+1$ sont dans $\M$, car $\M$ est complet, et contredit donc l'hypothse $|\M| = N$.  Ainsi, nous pouvons Žnoncer que
$$  P^+(m)\leqslant p_{N}\qquad (m\in\M).\eqdef{P+m}$$
Ici et dans la suite, nous dŽsignons par $P^+(m)$  (resp. $P^-(n)$)  le plus grand  (resp. le plus petit)  facteur premier d'un entier naturel~$n$ avec la convention $P^+(1)=1$  (resp. $P^-(1)=\infty$).  
\par 
 Notons $\M^*:=\big\{ [m,n] : m\in\M,\,n\in\M\big\}$. La majoration \eqref{majgenSf}
permet d'Žcrire
$$S \big( \M;C^\omega g_1\big)
\leqslant \sum_{\ell\in \M^*}\Bigg(\sum_{\di{m\mid \ell}{m\in \M}}
C^{\omega(\ell/m)}g_1\Big({\ell\over m}\Big) \Bigg)^2.\eqdef{majSMC}$$
\par 
Soit $f$ la fonction multiplicative dŽfinie par 
$$f(p)=\normalbaselineskip=20pt\cases{\dsp{ 1\over \sqrt{p}-1} & si $p\leqslant 2By$,\cr
\dsp C\sqrt{{ \log_3 N\over  \log N\log_2N}} 
\log \bigg({p\over B{\pi(y)}\log p}\bigg)   & si $p> 2By$.
\cr} $$
\par 
L'inŽgalitŽ de Cauchy-Schwarz fournit
$$\Bigg(\sum_{\di{m\mid \ell}{m\in \M}}
C^{\omega(\ell/m)}g_1\Big({\ell\over m}\Big)\Bigg)^2
\leqslant  \sum_{\di{d\mid \ell}{d\in \M}} f\Big({\ell\over d}\Big)
\sum_{\di{t\mid \ell}{t\in \M}} {C^{2\omega(\ell/t)} \over   f(\ell/t )}g_1\Big({\ell\over t}\Big)^2 \cdot\eqdef{CS}$$
Nous avons 
$$\sum_{\di{d\mid \ell}{d\in \M}}f\Big({\ell\over d}\Big)\leqslant F(\ell):=\prod_{p\mid \ell} (1+f(p))   $$
et
$$\sum_{p\leqslant 2By }  f(p)\ll {\sqrt{y}\over \log y}\ll{\sqrt{\log N\over \log_2N}} ,\eqdef{majsumy} $$
alors que, ds que $m\in \M$,
$$\sum_{\di{p\mid m}{p>2By  }}  f(p)\leqslant C
\sqrt{{ \log_3 N\over \log N\log_2N}} \sum_{\di{p\mid m}{p>y }} \log \bigg({p\over B\pi(y)\log p}\bigg)
\leqslant  C\sqrt{{ \log_3 N\log N\over \log_2N}}\cdot$$
Ainsi,  
$$F(m)\leqslant \exp\Big(\sum_{p\mid m} f(p)\Big)
\leqslant  \L(N)^{C+o(1)}  \qquad (m\in \M).\eqdef{maj1+f(p)}$$

En reportant \eqref{CS} dans le membre de droite de \eqref{majSMC}, nous obtenons donc
$$\eqalign{S \big(\M;C^\omega {g_1}\big)&\leqslant  
\sum_{\ell \in \M^*}\sum_{\di{m \mid \ell }{m \in \M}} { C^{2\omega(\ell /m)} \over   f(\ell /m )}g_1\Big({\ell\over m}\Big)^2 F(\ell) 
\leqslant \sum_{m\in \M }F(m) \sum_{ P^+( d)\leqslant p_N} { C^{2\omega( d)} g_1( d)^2F( d)\over   f( d )} \cr
&\leqslant |\M|\L(N)^{C+o(1)} \sum_{ P^+( d)\leqslant p_N}    { C^{2\omega( d)}  g_1( d)^2F( d)\over   f( d )}\cr
&= |\M|\L(N)^{C+o(1)}\prod_{ p\leqslant p_N}   \Big(1+{C^2(1+f(p))\over  (\sqrt{p}-1)^2f(p)}\Big)
\cr&\leqslant |\M|\L(N)^{C+o(1)} \exp\Big(\sum_{  p \leqslant p_N}    {C^2(1+f(p))\over (\sqrt{p}-1)^2f(p)}\Big)
.\cr}$$
Comme
$$\sum_{{ p\leqslant p_N}}  {1\over (\sqrt{p}-1)^2 }\ll \log_2N$$
et 
$$\eqalign{\sum_{p>2By} {C^2\over(\sqrt{p}-1)^2f(p)}&\leqslant \sqrt{{\log N\log_2N\over  \log_3 N }}\sum_{p>2By} {C\over  (\sqrt{p}-1)^2\log (p/B\pi(y)\log p)}
\cr&=\{C+o(1)\} \sqrt{{\log N\log_3N\over  \log_2 N }}\cr}$$
nous obtenons bien  \eqref{ssfc}.\qed
\bigskip 
\goodbreak
 
\paradeux{ComplŽtion de la preuve du \ref{th}}\drefdeux{complth1}

Nous adaptons la mŽthode de Bondarenko et Seip dans \citer{BS15}. 
La premire Žtape  consiste ˆ rŽduire la majoration de $S(\M)$ au cas d'un ensemble divisoriellement clos de mme cardinal --- voir \citer{ABS15}. Ë cette fin, nous modifions les ŽlŽments de $\M$ en altŽrant successivement, pour chaque \np\ $p$, les puissances de $p$ qui les divisent.

ConsidŽrons donc un nombre premier $p$. Notant $\M_p^*:=\{m/p^{v_p(m)}:m\in\M\}$, nous pouvons Žcrire  
$$\M=\bigcup_{t\in \M_p^*}\{ tp^{\nu_\alpha(t)}: 0\leqslant \alpha\leqslant r(t)\},\eqdef{decomM}$$
o $\nu(t):=\{\nu_\alpha(t)\}_{\alpha= 0}^{r(t)}$ dŽsigne la suite strictement croissante  des valuations $p$-adiques des ŽlŽments de $\M$ dont $t$ est la partie premire ˆ $p$. Il suit
$$S(\M; g_{0})=\sum_{m,n\in\M_p^*} g_{0}\Big(\frac{[m,n]}{(m,n)}\Big) \sigma_p( m,n),\eqdef{itep}$$
avec $$\sigma_p(m,n):=\sum_{\di{ 0\leqslant \alpha\leqslant r(m) }{ 0\leqslant \beta\leqslant r(n) } }{1\over  p^{|\nu_\alpha(m)-\nu_\beta(n)|/2}}\cdot$$
\par 
Commenons par Žtablir la majoration 
$$\max_{\di{ r(m) =r}{  r(n) = s}}\sigma_p ( m,n)\leqslant \sigma_p^*( r,s )\eqdef{majsigma}$$
o l'on a posŽ
$$\sigma_p^*( r,s ):=\normalbaselineskip=28pt\cases{\dsp r+1+  {2 r \over  \sqrt{p}-1}& si $r=s\geqslant  0$,\cr
\dsp
 r+1+ { 2r+1 \over  \sqrt{p}-1 } & si $r\geqslant 0,\,s= r+1$, \cr\dsp
 r+1+ {2r+2 \over  \sqrt{p}-1 } & si  $r\geqslant 0,\,s \geqslant  r+2$.\cr}\eqdef{sigp}$$  
 Ë cette fin, nous notons que, dans les deux premiers cas, l'optimum est atteint lorsque les ensembles $\nu(m)$ et $\nu(n)$ sont constituŽs des plus petits entiers possibles, autrement dit $\nu(m)= \{0, 1,\ldots, r\}$ et $\nu(n)=\{0, 1,\ldots, s\}$. 
Dans le troisime cas, la situation est plus compliquŽe. La configuration
$\nu(m)=\{0, 1,\ldots, r\}$, $\nu(n)=\{ j,  j+1,\ldots, j+s\}$  avec $j=\lfloor  (s-r)/2\rfloor$ est un exemple fournissant  l'optimum. 
Nous distinguons les cas $s=r$, $s=r+1$, et  $s \geqslant  r+2$. 
\par 
 Si $s=r$, nous avons 
$$\sigma_p ( m,n)\leqslant r+1+  
 \sum_{1\leqslant k\leqslant  r}{ 2(r+1 -k)\over   p^{ k /2}}\leqslant r+1+ { 2r   \over  \sqrt{p}-1 } =\sigma_p^*( r,r ), $$ 
alors que, si $s=r+1$, 
$$\sigma_p ( m,n)\leqslant r+1+   { 2r+1 \over   \sqrt{p}}+
 \sum_{1\leqslant k\leqslant  r}{ 2 r+1 -2k\over   p^{(s-r+k)/2}}\leqslant r+1+ { 2r+1   \over  \sqrt{p}-1 }=\sigma_p^*( r,r+1 ).  $$ 
\par 
Lorsque $s\geqslant r+2$, nous observons   que, pour chaque $k\geqslant 1$, le nombre de couples $(\alpha,\beta)$ tels que $|\alpha-\beta|=k$ n'excde pas $2(r+1)$. Il s'ensuit que
$$\sigma_p ( m,n)\leqslant r+1+  \sum_{ k\geqslant 1}{ 2r+2 \over   p^{k/2}}= r+1+ { 2r+2   \over  \sqrt{p}-1 }=\sigma_p^*( r,s ).  $$
Cela complte la preuve de \eqref{majsigma}.
\par \goodbreak

 Posons
$$h(p^\alpha,p^\beta):=\normalbaselineskip=15pt\cases{1 &si $\alpha=\beta$,\cr
g_1(p)& si $|\alpha-\beta|= 1$,\cr
0 & si $|\alpha-\beta|= 2$ et $\min(\alpha,\beta)\geqslant 1$,\cr
g_1(p) & si $|\alpha-\beta|= 2$ et $\min(\alpha,\beta)=0$,\cr
0 &si $|\alpha-\beta|\geqslant 3$,\cr
}\eqdef{hpp}$$ 
de sorte que
$$h(p^\alpha, p^\beta) ={1\over (\sqrt{p}-1)^{\min(| \alpha -\beta|,1) }}\qquad \big(|\alpha-\beta|\leqslant 2-\min\{ 1,\alpha\beta\}\big).$$
DŽfinissons alors les quantitŽs 
$$\sigma_p^{+}(r,s):=\sum_{\di{0\leqslant \alpha\leqslant  r }{0\leqslant \beta\leqslant s } }{h(p^\alpha, p^\beta) },\eqdef{hs}$$
\par 
Dans l'optique d'itŽrer une majoration issue de \eqref{itep}, nous allons ˆ prŽsent montrer que
$$\sigma_p^*(r,s)\leqslant \sigma_p^{+}\big(r,s \big)\qquad (r\geqslant 0,s\geqslant 0).\eqdef{majrt2}$$
Par symŽtrie, nous pouvons supposer $s\geqslant r$.\par 
Nous avons
$$\sigma_p^+(r,s)=u+{v+w\over \sqrt{p}-1},$$
avec
$$u:=\sum_{0\leqslant \alpha\leqslant r}1=r+1,\quad v:=\sum_{\tri{0\leqslant \alpha\leqslant r}{0\leqslant \beta\leqslant s}{|\alpha-\beta|=1}}1,\quad w:=\sum_{\tri{0\leqslant \alpha\leqslant r}{0\leqslant \beta\leqslant s}{|\alpha-\beta|=2,\,\alpha\beta=0}}1.$$
Si $r=s$, alors  $v=2r$, et $$w=\cases{0&$(r=s\leqslant 1)$,\cr
2&$(r=s\geqslant 2)$.\cr}$$
Ainsi $2r/(\sqrt{p}-1)\leqslant \sigma_p^+(r,r)-r-1\leqslant (2r+2)/(\sqrt{p}-1)$ et l'inŽgalitŽ \eqref{majrt2} est bien satisfaite.\par  
\par 
Si $s=r+1$, alors $v=2r+1$, et $$w=\cases{0&$(r=0,\,s=1$),\cr1&$(r=1,s=2)$,\cr 2 &$(r\geqslant 2,\,s=r+1).$}$$
D'o $(2r+1)/(\sqrt{p}-1)\leqslant \sigma_p^+(r,r+1)-r-1\leqslant (2r+3)/(\sqrt{p}-1) $, et nous obtenons encore \eqref{majrt2}. 
\par 
Si enfin $s\geqslant r+2$, alors $v=2r+1$, et $$w=\cases{1&$(r\leqslant 1,\,s\geqslant r+2)$,\cr2&$(r\geqslant 2,\,s\geqslant r+2)$.\cr}$$
Il suit $(2r+2)/(\sqrt{p}-1)\leqslant \sigma_p^+(r,s)-r-1\leqslant (2r+3)/(\sqrt{p}-1)$, ce qui confirme bien~\eqref{majrt2}.
\par 
  
\par 
Posons ˆ prŽsent $$\M_p=\bigcup_{n\in \M_p^*}\{np^\alpha:0\leqslant \alpha\leqslant  r(n) \},\eqdef{defMp}$$ de sorte que $|\M|=|\M_p|$.
Il rŽsulte de \eqref{itep}, \eqref{majsigma}, \eqref{hs} et \eqref{majrt2} que 
$$S(\M;g_0)\ \leqslant  \sum_{m,n\in\M_p^*} g_{0}\Big(\frac{[m,n]}{(m,n)}\Big) \sigma_p^+\big(r(m),r(n)\big).$$ 

Par itŽration, nous obtenons un ensemble $\M'$ divisoriellement clos tel que $|\M'|=|\M|$~et 
$$S( \M;g)\leqslant S^{*}(\M';h)\eqdef{S''M'} 
$$
avec $$S^{*}(\M';h):=
\sum_{m,n\in\M'} h(m,n)$$ et  $h$ est la fonction multiplicative de deux variables dŽfinies par \eqref{hpp}.

\par 
La deuxime Žtape, Žgalement prŽsente dans \citer{BS15}, consiste ˆ montrer que l'on peut essentiellement  se ramener au cas   o les ŽlŽments de $\M'$ sont sans facteur carrŽ.
\par 
 Pour simplifier les notations, nous considŽrons dans la suite un ensemble $\M$ divisoriellement clos et entreprenons de majorer $S^{*}(\M;h)$ tel que dŽfini plus haut. 
Il rŽsulte en particulier de notre hypothse que $\max_{m\in\M}\omega(m)\leqslant (\log N)/\log 2$.
\par 
DŽsignons par $k(m)$  le noyau sans facteur carrŽ d'un entier $m$.
Notant $\M_0$ l'ensemble des ŽlŽments sans facteur carrŽ de $\M$ et $\M_k:=\{ m\in \M : k(m)=k\}$ $(k\in\M_0)$, nous avons la  partition
$$\M=\bigcup_{k\in \M_0}\M_k.$$
et donc
$|\M|=N= \sum_{k\in \M_0}|\M_k|.$
Semblablement,
$$S^{*}(\M;h)=\sum_{j,k\in \M_0}S_{j,k}\eqdef{majS''MSjk}$$
avec
$$S_{j,k}:=\sum_{m\in \M_j,\,n\in\M_k} h (m,n ) \qquad (j,k\in\M_0).$$ 
Dans la somme intŽrieure, effectuons la d\'{e}composition 
$m=jm_1m_2,$ $ n=kn_1n_2,$
avec, pour $d:=(j,k)$,  $h(m,n)\neq0$,  
$$ m_2n_2\mid d^\infty,\quad m_1\mid  j/d ,\quad n_1 \mid  k/d.$$
puisque $h(m,n)$ est nul si $[m,n]/(m,n)$ a des facteurs cubiques. 
Nous avons 
$$(m_1,m_2)=(n_1,n_2)=(m_1,n_1)=1,\qquad \ell:={[m,n]\over (m,n)}={[j,k]\over (j,k)}{[m_2,n_2]\over (m_2,n_2)}m_1n_1.$$  
Cela implique
$$\eqalign{ h(m,n)&=  g_1\Big(\frac{[j,k]}{(j,k)}\Big) g_1\Big(\frac{[m_2,n_2]}{(m_2,n_2)}\Big)  \cdot\cr}$$
On vŽrifie cette formule par multiplicativitŽ en se restreignant au cas $(m,n)=(p^\alpha,p^\beta)$. Ainsi, par exemple,  si $(\alpha,\beta)= (2,0)$, nous avons $j=p$, $m_1=p $, $k=n_1=m_2=n_2=1$, et donc $h(p^2,1) =g_1(p) .$
\par Notons respectivement
$\M(j,k;m_2)$   et $\M(k,j;n_2)$   les ensembles constituŽs de toutes les valeurs admissibles de $m_1$ et $n_1$. 
Nous retrouvons l'inŽgalitŽ utilisŽe par Bondarenko et Seip au lemme 5 de \citer{BS15} en observant que $$\sum_{\di{m_1\in\M(j,k;m_2) }{ n_1\in \M(k,j;n_2)}}1\leqslant \sum_{\di{m_1\mid j/d}{n_1\mid k/d}} 1\leqslant   {2}^{\omega([j,k]/d)}.$$
 \par 
Posons
$$s(j,k;m_2)=\sum_{ m_1\in\M(j,k;m_2) } \mu(m_1 )^2 .\eqdef{notations}$$
Nous avons 
\ $$\eqalign{\sum_{\di{m_1\in\M(j,k;m_2) }{ n_1\in \M(k,j;n_2)}} \mu(m_1n_1)^2 =s(j,k;m_2)s(k,j;n_2).\cr}$$  
De plus, l'inŽgalitŽ de Cauchy--Schwarz fournit 
$$s(j,k;m_2)^2\leqslant |\M(j,k;m_2)|\sum_{m_1\mid j/d}1\leqslant |\M(j,k;m_2)|2^{\omega(j/d)}.\eqdef{majSjk}$$
Il s'ensuit que
$$\eqalign{S_{j,k}
&\leq\sqrt{2}^{\omega([j, k]/d)} g_1\bigg(\frac{[j,k]}{(j,k)}\bigg)S'_{j,k} 
  ,\cr}$$
avec
$$
S'_{j,k}:= \sum_{\di{\M(j,k;m_2)\neq\emptyset}{ \M({k,j};n_2)\neq\emptyset}}  
g_1\Big(\frac{[m_2,n_2]}{(m_2,n_2)}\Big) 
\sqrt{|\M(j,k;m_2)|\cdot |\M(k,j;n_2)|}.
$$

Une nouvelle application de l'inŽgalitŽ de Cauchy--Schwarz permet alors d'Žcrire
$$\eqalign{(S'_{j,k})^2&\leqslant \sum_{ {\di{\M(j,k;m_2)\neq\emptyset}{ \M({k,j};n_2)\neq\emptyset}}}  
g_1\Big(\frac{[m_2,n_2]}{(m_2,n_2)}\Big)  
  |\M(j,k;m_2)|  
\sum_{{\di{\M(j,k;m_2)\neq\emptyset}{ \M({k,j};n_2)\neq\emptyset}}}   
g_1\Big(\frac{[m_2,n_2]}{(m_2,n_2)}\Big)  
 |\M(k,j;n_2)| .\cr}
$$
Comme $\max\{\omega(j),\omega(k)\}
\leqslant (\log N)/\log 2$, toutes les valeurs considŽrŽes de $j$ et $k$ satisfont
$$G:=\prod_{p\mid  (j,k)} \big(1+{2 g_1(p) }\big)\ll G_N:=\e^{5\sqrt{(\log N)/\log_2N}}.\eqdef{defGN}$$ 
\par 
Pour chaque $m_2$  fixŽ tel que $\M(j,k;m_2)\neq\emptyset$, nous avons
$$\eqalign{\sum_{{\M(k,j;n_2)\neq\emptyset}}  
g_1\Big(\frac{[m_2,n_2]}{(m_2,n_2)}\Big) &\leqslant \prod_{p\mid (j,k) } \big(1+{2 g_1(p) }\big) =G\ll G_N.\cr}$$ 
Puisque
$$ \sum_{{\M(j,k;m_2)\neq\emptyset}} 
  |\M(j,k;m_2)|=|\M_{j}|,\qquad \sum_{{\M(k,j;n_2)\neq\emptyset}} 
  |\M(k,j;n_2)|=|\M_{k}|,  $$
nous obtenons
$$\eqalign{  S'_{j,k}
 & \ll   G_N\sqrt{|\M_{j}|\cdot|\M_{k}|}   \cr}
$$ 
et donc, en vertu de \eqref{majSjk},
$$S_{j,k}
 \leqslant G_N\sqrt{2}^{\omega([j,k]/(j ,k))}g_1\Big(\frac{[j,k]}{(j,k)}\Big)\sqrt{|\M_{j}|\cdot|\M_{k}|}.$$
 \par 
En reportant dans \eqref{majS''MSjk} et en appliquant le thŽorme 5 de \citer{ABS15}, nous obtenons
$$\eqalign{ S(\M)&\ll  G_N
\sum_{j,k\in \M}\mu(j)^2\mu(k)^2\sqrt{2}^{\omega([j,k]/(j,k))}g_1\Big(\frac{[j,k]}{(j,k)}\Big)\sqrt{|\M_{j}|\cdot|\M_{k}|}\cr&
\ll(\log N)G_N\sup_{\di{|\M| \leqslant N}{\M\subset \N_1}}S(\M,\sqrt{2}^{\omega}g_1).
\cr}$$ 
Le \ref{lemmessfc} fournit alors  la majoration du \ref{th}.
\qed
\medskip
\rem Cette dernire Žtape de la dŽmonstration vaut pour tout exposant $\alpha\in ]0,1[$: pour chaque ensemble $\M$ de cardinal $N$, nous avons 
$$S_\alpha(\M)\ll G_N(\alpha)^{O(1)}\sup_{\di{|\M'| = N}{\M'\subset \N_1}}S\big(\M',\sqrt{2}^{\omega}g_{2\alpha}\big).\eqdef{redgen}$$
avec
$$g_\alpha(n):={\mu(n)^2\over \prod_{p\mid n}( p^{\alpha/2}-1)},\qquad G_N(\alpha):=\exp\big\{ (\log N)^{1-\alpha}(\log_2N)^{-\alpha} \big\}.$$
Cela permet donc de majorer efficacement $S_{\alpha}(\M)$ lorsque $\dm<\alpha\leqslant \dm+o(1)$.
\smallskip
\bigskip
\goodbreak 
\paraun{Cas d'un ensemble de diviseurs : preuve du \ref{maj}}

Un calcul standard, effectuŽ dans \citer{LR17} lorsque $\alpha=1$, fournit
$$S_\alpha(\T_D)=\tau(D)\prod_{p^{\mu_p}\parallel D} \Big(1+{2\mu_p\over (1+\mu_p)p^\alpha}\sum_{0\leqslant k<\mu_p} {1-k/\mu_p\over p^{k\alpha}}
\Big).\eqdef{SaMD}
$$
Il suit, pour $\alpha:=1/2$,  
$$S(\T_D)\leqslant \tau(D)\exp\bigg\{\sum_{p^{\mu_p}\parallel D}  {2\mu_p\over (1+\mu_p)(\sqrt{p}-1)}   
 \bigg\} .\eqdef{majS2}$$
\par 
Commenons par traiter le cas des estimations de $\Gamma''(N)$.\par 
Lorsque $D$ est sans facteur carrŽ, nous dŽduisons de \eqref{SaMD} que 
$$
\exp\bigg\{\sum_{p\mid  D}  {1 \over  \sqrt{p} } 
 \bigg\}\prod_{p\mid D}\Big(1+{1\over 2p}\Big)^{-1}\leqslant 
{S(\T_D)\over \tau(D)}\leqslant  \exp\bigg\{\sum_{p\mid  D}  {1 \over   \sqrt{p} } 
 \bigg\}.\eqdef{majS1/2D}$$

\par 
Nous observons Žgalement que la contribution ˆ $S(\T_D) $ de chaque entier $m\in \T_D$   est alors une fonction de  $D$ seul : en effet, 
pour  $m\mid D$, $\mu(D)^2=1$ et $\alpha>0$, nous avons 
$$ \sum_{d \mid D}  {(d,m)^{\alpha}\over [d,m]^{\alpha}}=
\sum_{t \mid m }{ t^{ \alpha}\over m^{\alpha}}\sum_{s\mid D/m } {1 \over  s^{\alpha}} 
=\prod_{p \mid D } \Big(1+{1\over  p^{\alpha}}\Big)
.$$

\par 
La valeur maximale du membre de droite de \eqref{majS1/2D} lorsque  $\omega(D)=k$ est atteinte  lorsque  $D=D_y:=\prod_{p\leqslant y} p$ avec $\pi(y)=k=\fl{(\log N)/\log 2}.$
Gr‰ce ˆ l'estimation
$$\sum_{p\leqslant y}  {1 \over   \sqrt{p} } ={2\sqrt{y}\over \log y}\Big\{1+O\Big({1\over \log y}\Big)\Big\}={2\over \sqrt{\log 2}}\sqrt{\log N\over \log_2N}\Big\{1+O\Big({\log_3N\over \log_2N}\Big)\Big\},\eqdef{majsumpy}$$
nous obtenons 
$$\log \Gamma''_{1/2}(N)\leqslant  
{2\over \sqrt{\log 2}}\sqrt{\log N\over \log_2N}\Big\{1+O\Big({\log_3N\over \log_2N}\Big)\Big\} .$$
La minoration correspondante peut tre  obtenue de manire analogue. Nous omettons les dŽtails. 
\par 
ConsidŽrons ˆ prŽsent le cas de $\Gamma'(N)$.\par 
Pour optimiser la somme, disons $s(D)$, apparaissant au membre de droite de   \eqref{majS2} sous la contrainte $  \log \tau(D)\leqslant \log N$, nous Žcrivons la dŽcomposition canonique de $D$ sous la forme $D:=D_y=\prod_{p\leqslant y} p^{\mu_p}$.  La suite $\mu_p$ doit tre  dŽcroissante au sens large. En effet, comme la fonction $t\mapsto 2t/(t+1)$ est  croissante, si l'on avait $\mu_p>\mu_q$ avec $p<q$, en intervertissant les valeurs de $\mu_p$ et $\mu_q$, nous ne changerions pas la valeur de $\tau(D)$ tout en  accroissant $s(D)$. 
 
DŽfinissons alors une suite croissante au sens large $\{r_k\}_{k=1}^{\infty}\in[1,\infty[^\N$ et un entier $K\geqslant 1$ par les conditions
$$\mu_p=k\Leftrightarrow p\in I_k:=]y/r_{k+1},y/r_{k} ] \qquad \big(k\geqslant 1\big),\qquad K:=\inf\{k\geqslant 1:r_k>(\log y)^2\}.$$
Les intervalles $I_k$ peuvent tre vides: si la valeur $k$ n'est pas atteinte par les $\mu_p$, nous choisissons $r_{k+1}=r_k$. 
L'inŽgalitŽ \eqref{majS2} implique
$$\eqalign{\log&\Big({S(\M_{D})\over \tau(D)}\Big)\leqslant  \sum_{1\leqslant k\leqslant K} {2k\over 1+k } \sum_{y/r_{k+1}<p\leqslant y/r_{k} } {1\over \sqrt{p}-1}+ \sum_{ p\leqslant y/(\log y)^2 } {2\over \sqrt{p}-1}
\cr&\leqslant {4C_1\sqrt{y}\over \log y} \Big\{1+O\Big({1\over \log y}\Big)\Big\} ,\cr}
  $$
o la sŽrie $$C_1:=\sum_{k\geqslant 1} { k\over (1+k) }\bigg({1\over \sqrt{r_{k}}}-{1\over \sqrt{r_{k+1}}} \bigg)=
\sum_{k\geqslant 1} {1\over k(1+k) \sqrt{r_k}}   $$ est convergente.  
La condition $\tau(D)\leqslant N$ implique
$${y \over \log y} \Big\{C_2(K)+O\Big({1\over \log y}\Big)\Big\} \leqslant \log N \eqdef{C2}$$
o l'on a posŽ $$C_2(K):=  \sum_{1\leqslant k\leqslant K}\log (k+1)\bigg({1\over r_{k}}-{1\over r_{k+1}} \bigg)
=\sum_{1\leqslant k\leqslant K}{\log (1+1/k)\over r_k}-{\log (K+1)\over r_K}\cdot\eqdef{defC2K}$$ 
\par 
Pour Žvaluer le maximum de $C_1$ sous la contrainte \eqref{C2}, nous utilisons la mŽthode des  multiplicateurs de Lagrange  en considŽrant $\bfr:=\{r_k\}_{k=1}^{K}$ comme  variable et en nŽgligeant initialement les termes d'erreur ainsi que le deuxime terme du membre de droite de \eqref{defC2K}, qui est nŽgatif ou nul. Posons
$$
s(\bfr):=\sum_{1\leqslant k\leqslant K} {1\over k(1+k) \sqrt{r_k}},\qquad \psi(\bfr):=\sum_{1\leqslant k\leqslant K}{\log (1+1/k)\over r_{k}}-{\log N\log y\over y},$$ 
de sorte que nous cherchons le maximum de $s(\bfr)$ sous la contrainte $\psi(\bfr)\leqslant 0$. 
 Le   multiplicateur  de Lagrange $\lambda$ est donc dŽfini par les Žquations
$${\partial\over \partial \lambda}(s(\bfr)-\lambda \Psi(\bfr)\big)=0,\qquad{\partial\over \partial r_k}(s(\bfr)-\lambda \Psi(\bfr)\big)=0\qquad (1\leqslant k\leqslant K). $$
Ainsi  $\Psi(\bfr)=0$ et
$$- { 1\over 2 r_k^{3/2}k ( k+1 )  }  +{\lambda \over r_k^2 }\log (1+1/k)=0\qquad (1\leqslant k\leqslant K), $$
ce qui fournit 
$$r_k=\{2\lambda k(k+1)\log (1+1/k)\}^2.\eqdef{rk}$$
Nous choisissons alors  $\lambda=1/(4\log 2)$ de sorte que $r_1=1$.
Il s'ensuit que $r_k\asymp k^2$, $K\asymp \log y$ et 
$C_2(K)=C_2+O\big( (\log_2y)/(\log y)^{2}\big)$ avec $C_2:=\lim_{K\to +\infty } C_2(K).$\par 
De plus, compte-tenu de \eqref{C2}, nous avons  $y=(\log N)\big\{\log_2N+O(1)\}/C_2
$. Il suit
$$\eqalign{\log \Gamma'(N)\leqslant \log  \Big\{{S(\T_D)\over \tau(D)}\Big\}   \leqslant  {4C_1\over \sqrt{C_2}}\Big\{1+O\Big({\log_3N\over \log_2N}\Big)\Big\}\sqrt{\log N\over \log_2N}\cdot\cr}
  $$ 
Comme \eqref{rk} implique $4C_1/\sqrt{C_2}=B$,  nous obtenons bien la majoration contenue dans~\eqref{G'}. 
  \par 
Le mme choix des paramtres fournit alors la minoration correspondante.
\qed

\medskip 
\bigskip
\paraun{Valeurs extrmes de la fonction zta de Riemann}\drefun{vZ}
\paradeux{MŽthode de Bondarenko--Seip et sous-sommes de G‡l}
Dans leurs travaux sur le sujet, Bondarenko et Seip minorent les valeurs extrmes localisŽes de la fonction zta de Riemann --- et donc en particulier la quantitŽ \eqref{defZT} --- en Žvaluant les normes de formes quadratiques associŽes ˆ certaines sous-sommes de G‡l, soit
$$\gotS_\alpha( \M):=\sum_{\di{m,n\in\M}{n\mid m}} \Big({n\over m}\Big)^\alpha\qquad (\alpha>0).\eqdef{defSigmaa}$$
Nous avons bien, pour tout ensemble fini $\M$ et tout $\alpha>0$,
$$\gotS_\alpha( \M)\leqslant S_\alpha( \M),$$
puisque $[m,n]/(m,n)=m/n$ lorsque $n\mid  m$.
\par \goodbreak
Les minorations  issues, comme dans \citer{BS15}, de cette approche sont obtenues en construisant des ensembles $\M$ adŽquats qui sont divisoriellement clos, une condition naturelle pour obtenir de grandes valeurs. Or, il  est facile de constater que les sous-sommes de G‡l $\gotS ( \M)= \gotS_{1/2}( \M)$ sont dans ce cas d'un ordre de grandeur trs infŽrieur ˆ celui des sommes de G‡l maximales $S(\M)$.  
\Propl{sousGal}{Pour tout ensemble d'entiers $\M$ divisoriellement clos et de cardinal $N$, nous avons
$$\gotS (\M)  \leqslant  N \exp \Bigg\{ {2+o(1)\over \sqrt{\log 2}}\sqrt{\log N\over \log_2N}\Bigg\}\qquad (N\to\infty).$$ }
\nid
Pour tout $m$, nous avons 
$$\sum_{ n\mid m}  g_0\Big({m\over n}\Big)\leqslant \prod_{p\mid m} {1\over 1-1/\sqrt{p}}\cdot$$
De plus, si $\M$ de cardinal $N$ est divisoriellement clos nous avons vu que $\omega(m)\leqslant (\log N)/\log 2$ pour tout $m\in \M$, de sorte  que
$$\prod_{p\mid m} {1\over 1-1/\sqrt{p}}\leqslant \exp \Bigg\{ {2+o(1)\over \sqrt{\log 2}}\sqrt{\log N\over \log_2N}\Bigg\}.$$
La majoration annoncŽe en dŽcoule immŽdiatement.\qed
\medskip
 Posons 
$$\gotQ(\M):=\sup_{\di{\bfc\in \CC^N}{\no{\bfc}_2=1}}\abs{\sum_{\di{m,n\in\M}{n\mid m}}c_m\ov{c_n} \sqrt{{n\over m }}}.$$
Bondarenko et Seip Žtablissent dans \citer{BS15} l'existence d'un ensemble $\M$ divisoriellement clos et de cardinal $\leqslant N$ tel que
$$\gotQ(\M)\gg \L(N)^{ 1+o(1)}.\eqdef{minQM}$$
Cela met en Žvidence que la propriŽtŽ remarquable, mentionnŽe dans l'introduction, relative aux vecteurs propres des formes quadratiques associŽes aux sommes de G‡l  n'est plus valable pour les sous-sommes de type~$\gotS(\M)$, du moins lorsque $\M$ est divisoriellement clos. 
\par 
Nous avons vu plus haut qu'en vertu des lemmes 1 et 2 de \citer{BS15}, le supremum $\Gamma^*(N)$ est obtenu pour un ensemble strict. Nous savons  Žgalement  que tout ensemble strict  vŽrifie
$$\sup_{m\in\M}P^+(m)\leqslant p_N\qquad (|\M|=N).\eqdef{majP+bis}$$
  Notons $\A$ l'ensemble de toutes les parties finies de $\N_1$ satisfaisant \eqref{majP+bis}. Ainsi $\A$ contient ˆ la fois les ensembles considŽrŽs dans \citer{BS17} et ceux qui ont permis d'Žtablir la minoration du \ref{th}.\par 
Le r\'{e}sultat suivant, relatif aux ŽlŽments de $\A$, laisse supposer que la taille maximale de $\gotQ(\M)$ lorsque $\M$ dŽcrit les sous-ensemble de $\N_1$ de taille $N$ vaut $\Gamma^*(N)^{1/2+o(1)}$. Cette question demeure ouverte en l'Žtat.  Rappelons la dŽfinition des normes $Q(\M)$ en \eqref{defQM}. 
\par 
\Propl{quadra}{Soit $\M\in\A$, $|\M|=N$. Alors
 $$\gotQ(\M)^2
\ll (\log N)Q(\M).\eqdef{bof}$$
En particulier, 
$$\sup_{\di{\M\in\A}{|\M|=N}}\gotQ(\M)
=  \L(N)^{ 1+o(1)}.\eqdef{evalQM}$$}
\nid Soit $\bfc\in \CC^N$ tel que $\no{\bfc}_2=1$.
Posant
$$s_d:=\sum_{\di{m \in\M}{d\mid m}} c_mg_0\Big({m\over d}\Big),$$
de sorte que la quantitŽ ˆ majorer s'Žcrit $|\Lambda|^2$ avec $\Lambda:=\sum_{d\in\M}c_d\ov{s_d},$ nous appliquons 
l'inŽgalitŽ de Cauchy-Schwarz pour obtenir 
$$\eqalign{|\Lambda|^2&\leqslant  \sum_{d\in \M} |s_d|^2
 =\sum_{d\in \M}\sum_{\di{m,n\in \M}{d\mid m, d\mid n}}c_m\ov{c_n}{g_0(m)g_0(n)\over g_0(d)^2}
\leqslant   \sum_{m,n\in \M}g_0(m)g_0(n)|c_m {c_n}|\sum_{d\mid (m,n)}{1\over g_0(d)^2}\cdot\cr}$$
Or, pour tout entier $D$ tel que $P^+(D)\leqslant p_N$, nous avons 
$$\sum_{d\mid D}{1\over g_0(d)^2}\leqslant {1\over g_0(D)^2}\prod_{p^\nu\parallel D}{1-1/p^{-\nu-1}\over  1-1/p}\ll {\log N\over g_0(D)^2}\cdot$$
\par 
Il suit
$$\eqalign{|\Lambda|^2&\ll  (\log N)
\sum_{m,n\in \M}g_0\Big({[m,n]\over (m,n)}\Big)|c_m {c_n}|   \leqslant (\log N)  \sup_{\di{\bfc\in \CC^N}{\no{\bfc}_2=1}}\abs{\sum_{ {n,m\in \M} }c_m\ov{c_{n}} \sqrt{{(m,n)\over [m,n]}}} .\cr}\eqdef{majLambda}
$$
 Cela Žtablit bien \eqref{bof}.\par 
La variante du  thŽorme 5 de \citer{ABS15}\note{Voir \eqref{encth5}.} relative aux sous-ensembles de $\N_1$ permet, gr‰ce ˆ \eqref{maj^*Gamma}, de majorer le membre de droite de \eqref{majLambda} par $\L(N)^{2+o(1)}$. La minoration contenue dans \eqref{evalQM} a ŽtŽ Žtablie dans~\citer{BS17}.  
\qed

\medskip
\paradeux{Preuve du \ref{thzeta}}
Nous ferons appel ˆ la variante suivante du lemme 1 de \citer{BS17c}. Signalons incidemment qu'une erreur de signe s'est glissŽe dans ce dernier ŽnoncŽ  :  les deux membres du terme de droite de la formule~(6) de \citer{BS17c} doivent tre sŽparŽs par un signe $-$ et non $+$.
 \par 
Nous notons $z$ un nombre complexe gŽnŽrique et dŽfinissons implicitement ses parties rŽelle et imaginaire par $z=x+iy$. Lorsque $F\in L^1(\r)$, nous dŽfinissons la transformŽe de Fourier par la formule
$$\widehat F(\xi):=\int_\r F(x) \e^{-ix\xi} \d x\qquad (\xi\in\r).$$
\Propl{lemmaF}{Soient $\sigma\in ]-\infty,1[$ et $F$ une fonction holomorphe dans la bande horizontale $y=\im  z\in [\sigma-2,0]$, satisfaisant ˆ la condition de croissance
$$\sup_{ \sigma-2\leqslant y\leqslant 0}\big|F(z)\big|\ll{1\over x^2+1}\cdot\eqdef{condF}$$
Alors, pour tout $s=\sigma+it\in\CC$, $t\neq0$, nous avons
$$\eqalign{\int_\r \zeta(s+iu)&\ov{\zeta(s-iu)}F(u)\d u\cr
&=\sum_{k,\ell\geqslant 1}{ \widehat F(\log k\ell)\over k^{s} \ell^{\ov s}}-2\pi \zeta(1-2it) F(is-i)-2\pi \zeta(1+2it) F(i\ov s-i).\cr}$$}
\nid La formule rŽsulte d'une application standard du thŽorme des rŽsidus ˆ $w\mapsto\zeta(s+w)\zeta(\ov s+w)F(-iw)$, suivie d'un passage ˆ la limite. Nous omettons les dŽtails.
\qed
\medskip
\noi{\it DŽmonstration du \ref{thzeta}.} Soient $T>1$, $\kappa=1-\beta$, et $N:=\fl{T^{\kappa}}.$
ƒtant donnŽ un ensemble d'entiers  $\M$ de cardinal $N$, 
posons
$$\M_j:=\M\cap \ \big](1+1/T)^{j},(1+1/T)^{j+1}\big]\qquad (j\geqslant 0),\eqdef{defMj}$$ 
 notons $h_j:=\min\M_j$ lorsque $\M_j\neq\emptyset$, et dŽsignons par $\HH$ l'ensemble des $h_j$. DŽfinissons alors $r:\HH\to\r^+$ par la formule $r(h_j)^2:=\sum_{m\in\M_j}1$. Le  facteur de rŽsonance associŽ au problme  est choisi sous la forme 
$|R(t)|^2$ avec
$$R(t):=\sum_{h\in \HH}{r(h)\over h^{it}} \cdot$$
\par 
Nous avons trivialement
$$R(0)^2\leqslant   N \sum_{h\in \HH} r(h)^2\leqslant N  |\M| \leqslant N^2.\eqdef{inegR0}$$
 \par 
ConsidŽrons alors la densitŽ de Gauss $\Phi(t):=\e^{-t^2/2}$, de sorte que
$\widehat \Phi(\xi)=\sqrt{2\pi} \Phi(\xi).$
La mŽthode de rŽsonance repose classiquement \citer{BS17} sur une minoration du rapport $M_2(T)/M_1(T)$ o l'on a posŽ
$$\eqalign{M_1(T)&:=\int_\r |R(t)|^2\Phi\Big({t \over T}\Big)\d t,\quad
M_2(T):=\int_{T^{\beta}}^{T} |\zeta(\dm+it)|^2 |R(t)|^2\Phi\Big({t \over T}\Big)\d t.\cr}$$
 \par 
La technique de \citer{BS17c}, consistant,  une fois adapt\'{e}e ˆ notre situation,  ˆ remplacer, dans $M_2(T)$, le terme $|\zeta(\dm+it)|^2$ par un produit de convolution, permet un gain de prŽcision lorsque $\beta\leqslant \dm$.  
D'aprs le lemme 5 de \citer{BS17c}, nous avons 
$$\int_\r |R(t)|^2\Phi\Big({t\log T \over T}\Big)\d t \ll {T|\M|\over \log T} \cdot\eqdef{majavecM2}$$

Soit $\varepsilon\in\,]0,1[$. Posons
$$K(u):={\sin^2(\varepsilon u\log T)\over  \pi u^2\varepsilon \log T}\qquad (u\in \r)$$ de sorte que
$$\widehat K(\xi)=\Big\{1-{|\xi|\over 2\varepsilon \log T}\Big\}^+.$$
 Notons Žgalement, ˆ fins d'utilisation ultŽrieure, que $K$ est prolongeable en une fonction entire et satisfait \eqref{condF}. \par 
ConsidŽrons
$$\eqalign{\gotZ(t,u)&:=\zeta(\dm+it+iu)\zeta(\dm-it+iu)K(u),\cr
 I(T)&:=\int_\r|R(t)|^2\Phi\Big({t\log T \over T}\Big)\int_\r\gotZ(t,u) \d u\d t.\cr}
$$
Nous allons montrer dans un premier temps que les contributions ˆ $I(T)$ des domaines $|t|\leqslant 2T^\beta$ et $\{|t|>T/2 ,\,|u|>|t|/2\}$ sont nŽgligeables. Il s'ensuivra que seules des valeurs de 
$\zeta(\dm+iv)$ avec $v\in [ T^\beta,T]$ interviendront dans la contribution principale, ce qui permettra de majorer $I(T)$   en fonction de $Z_\beta(T)$. Dans un second temps, nous relierons $I(T)$ aux sommes de G‡l.
\par 
La majoration $$|\zeta(\dm +iv)|\ll (1+| v|)^{1/6}\qquad (v\in \r)\eqdef{majzeta}$$ fournit
$$\eqalign{\int_{|t|\leqslant 2T^\beta}&
\int_{|u|> T^\beta}  \gotZ(t,u)\d u\d t\ll
\int_{|t|\leqslant 2T^\beta}
\int_{|u|> T^\beta}   {|\zeta(\dm+it+iu)|\over 
|t|+|u|}{|\zeta(\dm-it+iu)|\over 
|t|+|u|} \d u\d t\ll T^{\beta}.\cr}$$
En appliquant l'inŽgalitŽ $2|ab|\leqslant |a|^2+|b|^2$ et en majorant classiquement le moment d'ordre $2$ sur $[-H,H]$ de $|\zeta(\dm+it)|$ par $\ll H\log H$, il vient
$$\eqalign{\int_{|t|\leqslant 2T^\beta}
\int_\r \gotZ(t,u)\d u\d t
&\ll T^\beta+\int_{|t|\leqslant 2T^\beta}\int_{|u|\leqslant T^\beta}
| \zeta(\dm+it+iu)|^2 K(u)\d u\d t \cr&
\ll T^\beta+\int_{|t|\leqslant 3T^\beta} 
| \zeta(\dm+it ) |^2  \d t\ll T^\beta\log T.\cr}$$
Il suit,  trivialement,
$$\eqalign{\int_{|t|\leqslant 2T^\beta} 
 |R(t)|^2\Phi\Big({t\log T \over T}\Big) \int_\r\gotZ(t,u)\d u\d t \ll R(0)^2 T^\beta\log T\ll
 |\M|  T^{\beta+\kappa}\log T.\cr}$$
\par 
De plus, la dŽcroissance rapide de $\Phi$ permet d'Žcrire
$$\int_{|t|>T/2} |R(t)|^2\Phi\Big({t\log T \over T}\Big)
\int_\r  \gotZ(t,u)\d u\d t\ll  R(0)^2.$$
Nous obtenons donc
$$\eqalign{I(T)+O\big(|\M| T \log T\big)=\int_{2T^{\beta}\leqslant |t|\leqslant T/2}|R(t)|^2\Phi\Big({t\log T \over T}\Big)\int_\r\gotZ(t,u)\d u \d t.
\cr}\eqdef{intdble}
$$
Comme, en vertu de \eqref{majzeta}, nous avons
$\int_{|u|>|t|/2}  \gotZ(t,u)\d u\ll 1,$ la contribution du domaine
$|u|>|t|/2$ au membre de droite de \eqref{intdble}  peut tre englobŽe par le terme d'erreur.
\par 
Ainsi l'intŽgrale intŽrieure peut tre restreinte au domaine  $T^\beta\leqslant |t\pm u|\leqslant T $.
Compte tenu de \eqref{majavecM2}, cela implique
$${T|\M|\over \log T}Z_\beta(T)^2
\gg I(T)+O\big(|\M| T \log T\big).\eqdef{majI}
$$
\par 
Ë ce stade, considŽrons la fonction
$$  G(z):=\sum_{ k,\ell\geqslant 1}{ \widehat K(\log k\ell)\over \sqrt{k\ell} (k/\ell)^{iz}}\cdot$$
L'application du \ref{lemmaF} avec $F=K$  fournit 
$$\eqalign{I(T)=I_1(T)+
I_2(T)+ I_3(T),\cr}\eqdef{est1IT}$$
o l'on a posŽ
$$\eqalign{I_1(T)&:= 
\int_\r  G(t)|R(t)|^2\Phi\Big({t\log T \over T}\Big)\d t,
\cr
I_2(T)&:=-2\pi \int_\r \zeta(1-2it) K(-t-\dm i)|R(t)|^2\Phi\Big({t\log T \over T}\Big)\d t,
\cr I_3(T)&:= -2\pi \int_\r \zeta(1+2it) K(t-\dm i)|R(t)|^2\Phi\Big({t\log T \over T}\Big)\d t.\cr} $$
En reportant dans \eqref{majI} l'estimation
$|I_2(T)|+|I_3(T)|\ll {T|\M|/\log T},$
qui rŽsulte des majorations \eqref{majzeta} et \eqref{majavecM2}, nous obtenons
$${T|\M|\over \log T} Z_\beta(T)^2
\gg I_1(T)+O\big(|\M|T\log T\big).\eqdef{prem}
$$
Cette estimation constitue la premire Žtape annoncŽe de la preuve.
\par \goodbreak
Afin de relier l'intŽgrale 
$I_1(T)$, aux sommes de G‡l, nous procŽdons comme dans \citer{BS17} --- \cf\ Žquation~(25).
Nous avons
$$\eqalign{I_1(T)&={T\sqrt{2\pi}\over \log T}\sum_{h,g\in \HH}r(h)r(g)\sum_{k,\ell\geqslant 1} 
{ \widehat K(\log k\ell)\over \sqrt{k\ell} }  \Phi\Big({T\over \log T}\log  {hk\over g\ell}\Big)
\cr&\gg  {T\over \log T}\sum_{1\leqslant k \ell\leqslant T^{\varepsilon }} 
{1\over \sqrt{k\ell} } \sum_{h,g\in \HH}r(h)r(g) \Phi\Big({T\over \log T}\log  {hk\over g\ell}\Big),\cr}$$
puisque le terme gŽnŽral de la sŽrie quadruple initiale est positif ou nul et $ \widehat K(\log k\ell)\geqslant 1/2$ si $k \ell\leqslant  T^{\varepsilon }$.
\par 
Selon la mŽthode dŽveloppŽe dans \citer{BS17}, nous observons que, lorsque $k$ et $\ell$ sont fixŽs et si $h\in\M_i$, $g\in\M_j$,   nous avons
$$\sum_{\di{m\in\M_i,\,n\in\M_j}{mk=n\ell}}1\leqslant \min\{ r(h)^2,r(g)^2\}\leqslant r(h)r(g),$$ 
et donc
$$\sum_{\di{m\in\M_i,\,n\in\M_j}{mk=n\ell}}\Phi\Big({T\over \log T}\log  {hn\over gm}\Big)\leqslant r(h)r(g)\Phi\Big({T\over \log T}\log  {hk\over g\ell}\Big).$$
Dans le membre de gauche, l'argument de $\Phi$ est $\ll 1/\log T$. Il suit, par sommation sur $i$ et $j$,
$$I_1(T)\gg{T\over \log T} \sum_{m,n\in \M}\sum_{\di{1\leqslant k \ell\leqslant T^{\varepsilon }}{mk=n\ell}} 
{1\over \sqrt{k\ell} }\cdot$$
En restreignant la somme intŽrieure aux couples $(k,\ell)$ tels que $$k=m/(n,m), \quad\ell=n/(n,m),\quad k\ell=[mn]/(m,n)\leqslant T^\varepsilon,$$ nous obtenons
$$\eqalign{I_1(T)&\gg {T\over \log T}\sum_{\di{m,n\in \M}{[m,n]/(m,n)\leqslant T^{\varepsilon}}} \sqrt{{(m,n)\over [m,n]}} \gg {T\over \log T}\Big(S_{1/2}(\M )-{S_{1/3}(\M)\over T^{\varepsilon/6}}\Big). \cr}$$

\par 
Reportons dans \eqref{prem}. Il vient
$$ Z_\beta(T)^2
\gg  {S_{1/2}(\M)\over |\M|}-{S_{1/3}(\M)\over |\M|T^{\varepsilon/6}} +O\big( (\log T)^2\big).\eqdef{minZ}
$$ 
Posant $y_\M:=\max_{m\in \M} P^{+}(m)$, nous avons pour tout $m\in \M$
$$\sum_{ n\in\M} {(m,n)^{1/3}\over [m,n]^{1/3}}\leqslant  
\prod_{p\leqslant y_\M} \Big(1+{2\over p^{1/3}-1}\Big)\ll  \exp\{ y_\M^{2/3}\}  $$ de sorte que 
$${S_{1/3}(\M)\over |\M|}\ll \exp  \big\{ y_\M^{2/3}\big\}.$$

 \par 
 SpŽcialisons alors $\M$ en choisissant l'ensemble dŽfini par \eqref{defM} avec $N=\fl{T^{1-\beta}}$, de sorte que 
$${S_{1/2}(\M)\over |\M|}\gg \L(T)^{2\sqrt{2(1-\beta)}+o(1)}.$$ Nous avons alors
$ y_\M\ll (\log T)^{6/5}$ ce qui implique que le terme nŽgatif de \eqref{minZ} est nŽgligeable. Il s'ensuit que
$$ Z_\beta(T)\gg \L(T)^{\sqrt{2(1-\beta)}+o(1)},$$
ce qui achve la dŽmonstration du \ref{thzeta}.\qed
\medskip\goodbreak

\rem   Un aspect essentiel de l'approche mise en Ïuvre dans la dŽmonstration prŽcŽdente repose sur l'introduction d'une intŽgrale relative ˆ un produit de deux facteurs zta, alors qu'un seul terme $|\zeta(\dm+it)|$ intervenait dans les travaux prŽcŽdents \citer{BS17}, \citer{BS17c}. Cela permet, via une majoration standard, de relier directement $Z_\beta(T)$ aux sommes de G‡l $S(\M)$ sans nŽcessiter l'introduction des formes quadratiques associŽes aux sous-sommes $\gotS(\M)$. Les coefficients $r(h)$ apparaissant  dans \citer{BS17} et \citer{BS17c} sont dŽfinis sous la forme $r(h_j)^2:=\sum_{m\in\M_j}f(m)^2$ 
o $f$ est une fonction multiplicative qui ne peut tre remplacŽe par $1$. L'introduction du carrŽ du module de zta permet 
de choisir  $f=\1$  et ainsi de minorer directement l'intŽgrale en fonction des  $S(\M)$,  sans appel aux normes $\gotQ(\M)$.
\par 
 \goodbreak

\paraun{Valeurs maximales de $|L(\dm,\chi)|$}

Nous nous proposons ici de dŽmontrer le \ref{thL}. Nous verrons ˆ cette occasion comment utiliser efficacement les minorations de sommes de G\'al dans ce problme. 
  Soit $\chi$ un caractre primitif modulo $q$. Alors  $\chi(-1) = (-1)^\nu$ o $\nu=\nu(\chi)\in\{ 0, 1\}. $ Ainsi $\chi$ est pair ou impair selon que $\nu(\chi)$ vaut 0 ou 1, et, posant
  $$W_\nu (x):={1\over 2\pi i} \int_{(2)}{\Gamma(\ft14+\dm s+\dm\nu)^2\over 
\Gamma(\ft14 +\dm\nu)^2sx^{s}} \d s, $$
o, ici et dans la suite nous notons $(\sigma)$ la droite verticale d'abscisse $\sigma$,
nous avons  
$$\vbs{L(\dm,\chi)}^2=2\sum_{k,\ell\geqslant 1}{\chi(k)\ov \chi (\ell)\over \sqrt{k \ell}}W_\nu \Big({\pi k\ell\over q}\Big)\qquad (\nu(\chi)=\nu).\eqdef{LW}$$
Ce rŽsultat classique  (\cf\ \citer{S07}, lemme 2)  rŽsulte de l'Žquation fonctionnelle satisfaite par la fonction $L(s,\chi)$ via un calcul standard de somme de Gauss associŽe ˆ un caractre primitif.
\par 
La formule dÕinversion de Mellin fournit l'expression suivante pour les fonctions $W_\nu$. \Propl{mellin}{Pour $x>0$ et $\nu\in\{0,1\}$, nous avons
$$ 
W_\nu (x)= 4\int_{x}^{\infty}
t^{\nu-1/2}\int_{ 0}^\infty {\e^{-v^2 -(t/v)^2}\over 
v\Gamma(\ft14 +\dm\nu)^2} \d v\d t
\qquad(\sigma>0).\eqdef{wnu}
$$
En particulier pour tout $x\geqslant 0$, nous avons $0\leqslant W_\nu (x)\leqslant 1$.}
\nid Nous avons classiquement (voir par exemple le corollaire II.0.14 
de \citer{Te08}) 
$$\e^{-t^2}={1\over 2\pi i} \int_{(\sigma)}\Gamma(s ){\dd s \over t^{2s}}\qquad(\sigma>0,\,t>0),$$
d'o, par changement de variables,
$$\e^{-t^2}t^{\nu+1/2}={1\over 4\pi i} \int_{(\sigma)}\Gamma(\dm s+\dm\nu +\ft14 ){\dd s \over t^s}$$
puis, en vertu de la formule de convolution adaptŽe ˆ la transformŽe de Mellin,
$$\eqalign{ {1\over 2\pi i} \int_{(\sigma)}\Gamma(\dm s+\dm\nu +\ft14 )^2{\dd s\over t^s }
&=4\int_{0}^\infty \e^{-v^2}v^{\nu+1/2}\e^{-(t/v)^2}(t/v)^{\nu+1/2}{\d v\over v}\cr&=4t^{\nu+1/2}\int_{0}^\infty \e^{-v^2 -(t/v)^2} {\d v\over v}\cdot
\cr}$$
Cela implique bien \eqref{wnu}, aprs division par $t$ et intŽgration sur $[x,\infty[$.\par
La positivitŽ de $W_\nu(x)$ est Žvidente ainsi que sa dŽcroissance par rapport ˆ $x$. D'aprs \citer{S07} --- formule (1.3a) --- nous avons $W_\nu(0)=1$ : cela implique immŽdiatement l'encadrement ŽnoncŽ. 
\qed
\goodbreak  
Rappelons ˆ prŽsent une variante, relative aux caractres primitifs, de la relation d'orthogonalitŽ des caractres de paritŽ donnŽe. D'aprs le lemme 1 de \citer{S07} (voir aussi le lemme 3 de \citer{BHB10}),  lorsque $(mn,q)=1$, nous avons 
$$\sumast_{\di{\chi\mod q}{\nu(\chi)=\nu}} \chi(m)\ov\chi(n)=\dm\sum_{d\mid (q,|m-n|)}\phi(d)\mu(q/d)+\dm(-1)^{\nu}\sum_{d\mid (q, m+n )}\phi(d)\mu(q/d),\eqdef{orthprim}$$
o, ici et dans la suite, l'astŽrisque indique que la sommation est restreinte aux caractres primitifs.
\par 
Supposons dorŽnavant que $q$ est un \np, de sorte que tout caractre non principal de  module $q$  est primitif.
Soit $\M=\M_q$ un ensemble  de cardinal  $N$ maximisant $S(\M)$ sous la contrainte que tous les ŽlŽments de $\M$ soient premiers ˆ $q$.
Telle que mise en Ïuvre par Soundararajan dans \citer{S08}, la mŽthode de rŽsonance consiste ˆ comparer les quantitŽs
$$\eqalign{V_1^+(q)&:=\sumast_{\chi \in X_q^+} |R_\chi|^2 ,\qquad\qquad
V_2^+(q):=\sumast_{\chi \in X_q^+} |R_\chi|^2 |L(\dm,\chi)|^2. \cr}$$
ƒtant donnŽ  un ensemble $\HH$ de reprŽsentants des classes de $\M$ modulo $q$, nous posons
 $$r(h)^2:=\sum_{\di{m\in\M}{m\md hq}} 1\qquad (h\in\HH),\qquad R_\chi:=\sum_{h\in \HH} r(h)\chi(h).\eqdef{defpoids}$$
L'inŽgalitŽ de Cauchy-Schwarz fournit, pour tout caractre $\chi$ de module $q$,
$$\vbs{R_\chi}^2\leqslant \vbs{R_{\chi_0}}^2\leqslant |\HH|\sum_{h\in \HH} r(h)^2\leqslant \min\{q-1, N \}N.\eqdef{inegRchi0}$$
La relation \eqref{LW} implique donc
$$V_2^+(q)= \sum_{j,h\in \HH}r(j)r(h)\sum_{\di{k,\ell\geqslant 1}{(k\ell,q)=1}}{1\over \sqrt{k\ell}}  W_0\Big({\pi k\ell\over q}\Big)\sigma_q,\eqdef{V2+}$$
avec, en vertu de \eqref{orthprim},
$$\sigma_q:=  2\sumast_{\chi \in X_q^+} \chi(jk)\ov \chi(h\ell)=  \sum_{d\mid (q,  jk+ h\ell )}\phi(d)\mu(q/d)+   \sum_{d\mid (q,| jk- h\ell|)}\phi(d)\mu(q/d)
   . $$ 
\par 
Comme $q$ est un \np, nous avons  
$$\sigma_q=\normalbaselineskip=15pt\cases{-2 & si $q\nmid (jk- h\ell)( jk+ h\ell),$\cr
 q-3 & si $q\nmid jk- h\ell$ et $q\mid  jk+ h\ell,$\cr
 q-3  & si $q\mid jk- h\ell$ et $q\nmid  jk+ h\ell,$\cr
  2(q-2) & si $q\mid jk- h\ell$ et $q\mid  jk+ h\ell.$\cr} $$
De plus, la relation $W_0(x)\ll 1/(x+1)^2$, valable pour $x\geqslant 0$, fournit
$$\sum_{\di{k,\ell\geqslant 1}{(k\ell,q)=1}}{1\over \sqrt{k\ell}}  W_0\Big({\pi k\ell\over q}\Big)\ll \sqrt{  q}\log  q.$$
Nous en dŽduisons 
$$\eqalign{V_2^+(q)+ 2\vbs{R_{\chi_0}}^2\sum_{\di{k,\ell\geqslant 1}{(k\ell,q)=1}}{1\over \sqrt{k\ell}}  W_0\Big({\pi k\ell\over q}\Big)& 
=\sum_{j,h\in \HH}r(j)r(h)\sum_{\di{k,\ell\geqslant 1}{(k\ell,q)=1}}{1\over \sqrt{k\ell}}  W_0\Big({\pi k\ell\over q}\Big)\big(\sigma_q+2\big)\cr
&\geqslant  (q-  1) \sum_{\di{k,\ell\geqslant 1}{(k\ell,  q)=1}}{1\over \sqrt{k\ell}}  W_0\Big({\pi k\ell\over q}\Big)\sum_{\di{j,h\in \HH}{q\,\mid\,  \abs{jk- h\ell}}}r(j)r(h),\cr}$$
o nous avons utilisŽ, de manire cruciale, la positivitŽ de $W_0$ Žtablie au \ref{mellin}. Il s'ensuit que
$$\eqalign{V_2^+(q)\geqslant (q- 1) \sum_{\di{k,\ell\geqslant 1}{(k\ell,  q)=1}}{1\over \sqrt{k\ell}}  W_0\Big({\pi k\ell\over q}\Big)\sum_{\di{j,h\in \HH}{q\,\mid\,  \abs{jk- h\ell}}}r(j)r(h)+O\Big(\min\{q, N \}N\sqrt{q}\log q\Big).\cr}$$     
   
 Lorsque $kj\md{\ell h}q$, nous avons
$$r(j)r(h)\geqslant \min\{r(j)^2,r(h)^2\}\geqslant \sum_{\di{m,n\in\M,\,km=\ell n}{m\md jq,\,n\md hq}}1,\eqdef{inegrjrh}$$  
de sorte que  
$$\eqalign{V_2^+(q)&\gg q\sum_{\di{k,\ell\leqslant \sqrt{q}}{(k\ell,q)=1}}{1\over \sqrt{k\ell}}
\sum_{\di{m, n\in\M }{km = \ell n }}1+O\Big(\min\{q, N \}N\sqrt{q}\log q\Big)\cr&\gg q
\sum_{\tri{m, n\in\M }{m/(m,n),n/(m,n)\leqslant \sqrt{q}}{(mn,q)=1}}\sqrt{(m,n)\over [m,n]}+O\Big(\min\{q, N \}N\sqrt{q}\log q\Big),\cr}$$
o nous avons utilisŽ la minoration $W_0(x)\gg 1$, valable pour $x\in [0,\pi]$.

Choisissons alors $N=\fl{\sqrt{q}}$. Nous avons construit dans la dŽmonstration du \ref{th} un ensemble $\M=\M_N$  de cardinal $N$ et tel que $S(\M_N)\geqslant N\L(N)^{2\sqrt{2}+o(1)}$. Or, tous les facteurs premiers $p$ des ŽlŽments de $\M_N$ vŽrifient $\log N<p\leq (\log N)^2\leq N$ pour $N$ suffisamment grand. Il s'ensuit que, si
$m\in \M_N$ et $N= \fl{\sqrt{q}}$, nous avons nŽcessairement $(m,q)=1$.
Il suit
$$\sum_{\tri{m, n\in\M_N }{m/(m,n),n/(m,n)\leqslant  \sqrt{q} }{(mn,q)=1}}\sqrt{(m,n)\over [m,n]}\geq S_{1/2}(\M_N)-{1\over N^{1/6}}S_{1/3}(\M_N)
\gg S_{1/2}(\M_N)\gg q^{1/2}\L(q)^{2+o(1)}$$
o l'avant-dernire minoration rŽsulte d'une manipulation semblable ˆ celle qui a permis de minorer le membre de droite de \eqref{minZ}.
 Pour ce choix de l'ensemble $\M$, nous avons donc 
$$V_2^+(q)\gg q^{3/2}\L(q)^{2+o(1)}.$$
De plus,
$$ V_1^+(q)\leqslant \sum_{\chi\,(\mod q)}|R_\chi|^2\leqslant q\sum_{m\in\M} r(m)^2\leqslant qN\leqslant q^{3/2}.   $$
Nous obtenons  donc la minoration 
$$L_q^+\geqslant \sqrt{V_2^+(q)/V_1^+(q)}\geqslant \L(q)^{1+o(1)}.$$
Cela achve la dŽmonstration du \ref{thL}.
\medskip

\bigskip

\paraun{Sommes de caractres}

\paradeux{Un lemme de localisation et  coprimalitŽ}
 
Le rŽsultat suivant est une consŽquence facile de la construction employŽe pour Žtablir le \ref{th}. 
 \Propc{corth}{Soient $q\geqslant 1$ et $N\geqslant 3$ tels que  $\pi(\log N\log_2N)\geqslant \omega(q)$. Il existe alors un ensemble d'entiers $\M_q$ vŽrifiant $$|\M_q|\leqslant N, \quad\max\M_q\leqslant 2\min\M_q,  \quad \big(\txt\prod_{m\in\M_q}m,q\big)=1,\quad S(\M_q)\geqslant |\M_q|\L(N)^{2\sqrt{2} +o(1)}.\eqdef{minSloc}$$} 
\nid  
ConsidŽrons l'ensemble $\M$  de cardinal $\leqslant N$ construit dans la dŽmonstration du \ref{th}. En faisant tendre  $\beta$ vers $2\sqrt{2}$ dans \eqref{minSMb}, nous avons donc 
$$S(\M)\geqslant |\M|\L(N)^{2\sqrt{2}+o(1)}.\eqdef{minSMc}$$
La premire Žtape consiste, sans changer $|\M|$ et sans diminuer $S(\M)$, ˆ altŽrer $\M$ de manire ˆ rendre tous ses ŽlŽments premiers ˆ $q$.  Ë cette fin, observons que $P^-(m)>\log N\log_2N$ pour tout $m$ de~$\M$,  dŽsignons par $p_j$ le $j$-ime \np \ ne divisant pas $q$, dŽcomposons $q=q_1q_2$ avec $P^+(q_1)\leq  \log N\log_2N$, $P^-(q_2)>  \log N\log_2N$, et posons $r=\omega(q_2)$. Ainsi $p_j\nmid \prod_{m\in\M}m$ $(1\leqslant j\leqslant r)$ alors que,
 par hypothse,  $ r=\omega(q_2)\leqslant \pi(\log N\log_2N)-\omega(q_1) $. DŽsignons par $\ell_1,\ldots, \ell_r $ les facteurs premiers de~$q_2$,  et associons ˆ chaque ŽlŽment $m$ de $\M$ l'entier $$m^\dagger:=m\prod_{1\leqslant j\leqslant r}\Big({p_j\over \ell_j}\Big)^{v_{\ell_j}(m)}$$ 
L'ensemble $\M^\dagger:=\{m^\dagger:m\in\M\}$ vŽrifie clairement $|\M^\dagger|=|\M|$. De plus, $S(\M^\dagger)\geq S(\M)
$ puisque $(m^\dagger,n^\dagger)/[m^\dagger,n^\dagger]\geqslant (m,n)/[m,n]$ pour tous $m,\,n$ de $\M$. 
\par 
 Dans une seconde Žtape, montrons que, sans diminuer significativement $S(\M^\dagger)/|\M^\dagger|$, on peut imposer que $\M^\dagger$ soit inclus dans un intervalle dyadique. Tout ŽlŽment $m$ de $\M^\dagger$ est  de la forme 
$m=(d / b)\prod_{1\leqslant k\leqslant (\log_2N)^\gamma}N_k^\dagger $ o l'on a posŽ $$N_k^\dagger:=\prod_{\di{p\in I_k}{p\nmid q}} p\prod_{\di{\ell_j\in I_k}{1\leqslant j\leqslant r} }p_j $$ et 
$$\log (d  b)\leqslant  \sum_{1\leqslant k\leqslant (\log_2N)^\gamma}
{au^{k+1}(\log N)^2\log_2N\over  k^2\log_3N}\ll (\log N)^3.$$
Nous pouvons donc scinder $\M^\dagger$ en au plus  $J\ll(\log N)^3$ sous-ensembles $\M_j$ inclus dans un intervalle dyadique $]M_j,2M_j]$. 
Or
$$\eqalign{S(\M^\dagger)&=\sum_{d} \phi(d) \Bigg(\sum_{\di{m\in \M^\dagger}{d\mid m}}{1\over \sqrt{m}}\Bigg)^2
\leqslant  J \sum_{d} \phi(d) \sum_{1\leqslant j\leqslant J}\Bigg(\sum_{\di{m\in \M_j}{d\mid m}}{1\over \sqrt{m}}\Bigg)^2
\cr&\leqslant  J\sum_{1\leqslant j\leqslant J}S(\M_j)\leqslant J^2\max_{1\leqslant j\leqslant J}S(\M_j).\cr}$$ 
Cela implique immŽdiatement le rŽsultat annoncŽ.
\qed
\medskip
\goodbreak

\paradeux{DŽmonstration du \ref{hough}}

 Soit $N$ un entier tel que $\pi(\log N\log_2N)\geqslant \omega(q)$.   Nous appliquons la mŽthode de rŽsonance en considŽrant
$$  
W_1 (q):=\sum_{\di{\chi\mod q}{\chi\neq \chi_0}} |R_\chi|^2 ,\qquad 
W_2 (q):=\sum_{\di{\chi\mod q}{\chi\neq \chi_0}} |R_\chi|^2 \vbs{S(x,\chi)}^2,  $$  
o $R_\chi$  est dŽfini comme en \eqref{defpoids} pour  l'ensemble $\M$ de cardinal $\leqslant N$ construit au \ref{corth}. 
\par 
L'orthogonalitŽ des caractres fournit d'abord
$$W_1(q)\leqslant \phi(q)|\M|.\eqdef{majW1}$$
Pour minorer $W_2(q)$, nous observons d'une part que, d'aprs \eqref{inegRchi0},
$$W_2 (q)+ |R_{\chi_0}|^2 \vbs{S(x,\chi_0)}^2\leqslant  W_2 (q)+O(|\M|^2x^2),\eqdef{majW2}$$
et, d'autre part, que l'orthogonalitŽ des caractres et l'inŽgalitŽ \eqref{inegrjrh} permettent d'Žcrire
$$\eqalign{
W_2 (q)+ |R_{\chi_0}|^2 \vbs{S(x,\chi_0)}^2
&=\phi(q)\sum_{\di{1\leqslant k,\ell\leqslant x}{(k\ell,q)=1}}\sum_{\di{j,h\in \HH}{q\,\mid\,  \abs{jk- h\ell}}}r(j)r(h)
\cr&\geqslant  \phi(q)\sum_{m,n\in \M}\sum_{\tri{1\leqslant k,\ell\leqslant x}{(k\ell,q)=1}{ mk= n\ell}}1. \cr}\eqdef{minW2}$$
\par 
Restreignons la somme intŽrieure, disons $\sigma(m,n)$ aux couples $(k,\ell)=(sn/(m,n),sm/(m,n))$ tels que $(s,q)=1$.
 Sous l'hypothse $\max\{ m,n\}\leqslant \sqrt{2mn}\leqslant  x(m,n)/(\log q)^A$, o $A$ est une constante absolue convenablement choisie, nous avons donc, en vertu, par exemple, du lemme 3.9 de \citer{BT05},
$$\sigma(m,n)\gg {\varphi(q)x (m,n)\over q\sqrt{mn}}= {\varphi(q)x\over q}\sqrt{{(m,n)\over  [m,n]}} .$$ 

Il rŽsulte alors de  \eqref{majW1}, \eqref{majW2} et \eqref{minW2} que
 $$\eqalign{\Delta(x,q)^2&\geqslant {W_2(q)\over W_1(q)}\gg  {\varphi(q)x\over q|\M|}\sum_{\di{m,n\in\M}{[m,n]/(m,n)\leqslant  x^2/2}} \sqrt{(m,n)  \over [m,n]}+O\Big({x^2|\M|\over {\varphi(q)}}\Big)\cr&\gg {\varphi(q)x\over q|\M|}\bigg\{S(\M)-{ 2S_{1/2-\eta}(\M)\over x^{2\eta}}\bigg\}+ O\Big({x^2N\over {\varphi(q)}}\Big),\cr} \eqdef{minD}$$ 
pour tout $\eta>0$.
\par 
Choisissons par exemple $\eta=\varepsilon/3$, o, par hypothse, $\varepsilon>0$ est tel que $\log x>(\log q)^{1/2+\varepsilon}$. 
Notant  ˆ nouveau $y_\M:=\max_{m\in \M} P(m)\leqslant (\log N)^{1+o(1)}$, nous avons, pour chaque $m$ de $ \M$,
$$\sum_{ n\in\M} \bigg({(m,n)\over [m,n]}\bigg)^{1/2-\eta}\leqslant  
\prod_{p\leqslant y_\M} \Big(1+{2\over p^{1/2-\eta}-1}\Big)\ll  \exp\{ y_\M^{1/2+\eta}\} \ll\exp\big\{(\log N)^{1/2+2\varepsilon/3}\big\} $$ de sorte que 
$${S_{1/2-\eta}(\M)\over |\M|}\ll \exp\big\{(\log N)^{1/2+2\varepsilon/3}\big\}.$$

 Pour le choix   $N:= K\fl{q/x}$, o $K$ est une constante absolue assez grande. En vertu de l'hypothse $x\leqslant q/\e^{(1+\varepsilon)\omega(q)}$, la condition $\pi(\log N\log_2N)>\omega(q)$ du \ref{corth} est bien satisfaite. Nous dŽduisons donc~\eqref{H+} de \eqref{minD} en y reportant la minoration \eqref{minSloc}. 

\bigskip
\par 
\noi{\bf Remerciements.} Le premier auteur  prend  plaisir ˆ remercier Marc Munsch pour d'intŽressantes discussions concernant les valeurs extrŽmales des fonctions $L$.

\bigskip\medskip\goodbreak
\centerline {\twelvebf Bibliographie}
\bigskip
{\eightpoint\leftskip9truemm

\bibtem{Ai16} C. Aistleitner, Lower bounds for the maximum of the Riemann zeta function along vertical lines, {\it Math. Ann. \bf365} (2016), \numeros 1-2, 473Ð496.\par 
\bibtem{AKMP18} C. Aistleitner, K. Mahatab, M. Munsch \& A. Peyrot,
On large values of $L(\sigma,\chi)$, {\it Quart. J. Math.}, ˆ para"tre.

\bibtem{ABS15} C. Aistleitner, I. Berkes \& K. Seip,   GCD sums from Poisson integrals and systems of dilated functions, {\it  J. Eur. Math. Soc.} {\bf  17} (2015), \numero6, 1517--1546. 

\bibtem{BR77}
R. Balasubramanian \& K. Ramachandra,   On the frequency of TitchmarshÕs phenomenon for $\zeta(s)$, III, {\it Proc. Indian Acad. Sci.} Sect. A {\bf 86} (1977), 341--351. 

\bibtem{BHS16} A. Bondarenko, T. Hilberdink \& K.  Seip,  G‡l-type GCD sums beyond the critical line, {\it  J. Number Theory}  {\bf 166} (2016), 93--104.

 \bibtem{BS15}
 A. Bondarenko,  K. Seip,  GCD sums and complete sets of square-free numbers, {\it  Bull. Lond. Math. Soc.} {\bf 47} (2015), \numero1, 29--41.

 \bibtem{BS17}
 A. Bondarenko,  K. Seip,  Large greatest common divisor sums and extreme values of the Riemann zeta function,  {\it  Duke Math. J.} {\bf 166} (2017), \numero9, 1685--1701.

\bibtem{BS17b}
 A. Bondarenko \&  K. Seip,   
    Note on the resonance method for the Riemann zeta function, in: ``Operator Theory: Advances and Applications"  261 (2018), 121--140, Birkh\"auser Verlag.

\bibtem{BS17c}
 A. Bondarenko \&  K. Seip, 
    Extreme values of the Riemann zeta function and its argument, {\it Math. Ann.} (2018), ˆ para"tre.

\bibtem{BT05} R. de la Bretche \& G. Tenenbaum, PropriŽtŽs statistiques des entiers friables, {\it Ramanujan J. \bf9} (2005),
139--202.\par 
\bibtem{BHB10} H.M. Bui \& D.R. Heath-Brown, A note on the fourth moment of Dirichlet L-functions, {\it Acta arith. \bf141}, \numero4 (2010), 335--344.
\par 
\bibtem{G49}	I.S. G‡l, A theorem concerning Diophantine approximations, {\it Nieuw Arch. Wiskunde} {\bf 23} (1949), 13--38. 

\bibtem{GS01}
A. Granville \& K. Soundararajan, Large character sums, {\it J. Amer. Math. Soc.} {\bf 14} (2) (2001), 365--397.

\bibtem{Hi09}	T. Hilberdink, An arithmetical mapping and applications to $\Omega$-results for the Riemann zeta function, {\it Acta Arith.} {\bf 139} (2009), 341--367.

\bibtem{H11}
B. Hough, The resonance method for large character sums, disponible sur https://arxiv.org/abs/1109.1786, (2011). 

\bibtem{H13}
B. Hough, 
The resonance method for large character sums, {\it  Mathematika} {\bf  59} (2013), no. 1, 87--118.

\bibtem{H16}
B. Hough,  The angle of large values of
$L$-functions. {\it J. Number Theory} {\bf  167} (2016), 353--393. \par

\bibtem{L11} Y. Lamzouri,   On the distribution of extreme values of zeta and $L$-functions in the strip $\dm<\sigma<1$. {\it Int. Math. Res. Not.} 2011, no. 23, 5449--5503.
\par 

\bibtem{LR17} M. Lewko \& M. Radziwi\l\l,   Refinements of G‡l's theorem and applications, {\it Adv. Math.} {\bf 305} (2017), 280Ð297.

\bibtem{M77} H. L. Montgomery, Extreme values of the Riemann zeta function, {\it Comment. Math. Helv.} {\bf 52} (1977), 511--518


\bibtem{SS16}
E. Saksman, K. Seip,  Some open questions in analysis for Dirichlet series, in: {\it Recent progress on operator theory and approximation in spaces of analytic functions}, 179--191, Contemp. Math.~{\bf 679}, Amer. Math. Soc., Providence, RI, 2016. 

\bibtem{S07} K. Soundararajan, 
The fourth moment of Dirichlet L-functions, in: Analytic Number Theory: A Tribute to Gauss and Dirichlet, Clay Math. Proc. 7, Amer. Math. Soc., Providence, RI, 2007, 239--246.

\bibtem{S08} K. Soundararajan, Extreme values of zeta and $L$-functions, {\it Math. Ann.} {\bf 342} (2008), 467Ð486. 

\bibtem{Te08} G. Tenenbaum, {\it Introduction ˆ la thŽorie analytique et probabiliste des
nombres}, quatrime Ždition, coll.~ƒchelles, Belin, 2015.\par
}

\bigskip
{\leftskip9mm\rightskip-2cm\sevenrm
\gutter=4cm \doublecolumns
 \obeylines \baselineskip=7pt
RŽgis de la Bretche
Institut de MathŽmatiques de Jussieu
    UMR 7586
UniversitŽ Paris Diderot-Paris 7
Sorbonne Paris CitŽ, 
Case 7012, F-75013 Paris
France
GŽrald Tenenbaum
 Institut ƒlie Cartan
 UniversitŽ de Lorraine
 BP 70239
 54506 Vand{\oe}uvre-ls-Nancy Cedex
 France
\phantom a
\singlecolumn
\par
}

\bye